\documentclass[12pt]{article}

\usepackage{bm, amsmath, amssymb}
\usepackage[dvips]{graphicx}

\newcommand{\eps}{\varepsilon}
\def\bx{\vec x}
\def\bu{\vec u}
\def\bv{\vec v}

\def\dt{\partial_t\/}
\def\dx{\partial_x\/}
\def\dy{\partial_y\/}
\def\NN{\mathbb N}
\def\ZZ{\mathbb Z}
\def\RR{\mathbb R}
\def\calf{\mathcal{F}}

\def\tr{\operatorname{Tr}}
\def\id{\operatorname{id}}
\def\Div{\operatorname{div}}
\newcommand\Ric{\operatorname{Ric}}
\newcommand\vol{\operatorname{vol}}

\newcommand{\eq}{\hspace*{-2mm}&=&\hspace*{-2mm}}
\newcommand{\plus}{\hspace*{-2mm}&+&\hspace*{-2mm}}

\newtheorem{cor}{Corollary}
\newtheorem{df}{Definition}
\newtheorem{ex}{Example}
\newtheorem{rem}{Remark}
\newtheorem{lem}{Lemma}
\newtheorem{prop}{Proposition}
\newtheorem{thm}{Theorem}
\newcounter{rembango}

\title{Extrinsic geometric flows\\ on foliated manifolds, I}

\author{Vladimir Rovenski$^{(1,2)}$ and
        Pawe\l\ Walczak$^{(2)}$\thanks{E-mail: rovenski@math.haifa.ac.il
                                           and pawelwal@math.uni.lodz.pl} \\
        (1) University of Haifa, (2) Uniwersytet \L\'odzki}

\date{}

\begin{document}

\maketitle

\begin{abstract}
We study deformations of Riemannian metrics on a given manifold equipped
with a codimension-one foliation subject to quantities expressed in
terms of its second fundamental form.
We prove the local existence and uniqueness theorem and estimate the existence time of solutions for some particular cases.
The key step of the solution procedure is to find (from a system of quasilinear PDEs) the principal curvatures of the foliation.
Examples for
extrinsic Newton transformation flow, extrinsic Ricci flow, and applications to foliations on surfaces are given.


\end{abstract}


\section*{Introduction}

Some of the most striking recent results in differential geometry and topology are related
to the \textit{Ricci flow} (see, among the others, \cite{to}),
that is the deformation $g_t$ of a given Riemannian metric $g_0$ on
a manifold $M$ subject to the partial differential equation (PDE)
\begin{equation*}
 \dt g_t = -2\Ric_t,
\end{equation*}
where $\Ric_t$ is the Ricci tensor on the Riemannian manifold $(M, g_t)$.
(Note that there is a -- due to Hamilton and based on suggestions of Yau
-- strategy of using the Ricci flow in proving
the famous Poincar\'e Conjecture and the Thurston Geometrization Conjecture
for 3-manifolds in a series of preprints by G. Perelman.)
 On the other hand, since at least 30 years there is continuous interest in the study
of \textit{Mean Curvature Flow} (MCF), i.e., the variation of immersions $F_0:\bar M\to M$ of manifolds $\bar M$
into Riemannian manifolds $(M, g)$ subject to the \textit{mean curvature},
i.e., to the PDE
\begin{equation*}
 \dt F_t = H_t,
\end{equation*}
where $H_t$ is the mean curvature vector of $F_t(\bar M)$
Also, the second author \cite{wa} considered the foliated version of the MCF: foliations of Riemannian
manifolds which are invariant under the flow of the mean curvature vector
of their leaves.
The recent years have seen a growing interest in \textit{Geometric Flows} $\dt g_t = h(b_t)$ of different types.

 In the paper we introduce \textit{Extrinsic Geometric Flow} (EGF) as
deformations of Riemannian metrics on a given manifold $M$ equipped with a codimension-one foliation subject to
quantities expressed in terms of its second fundamental form (suitably
extended to a symmetric tensor field of type $(0,2)$ on $TM$).
Authors propose EGF as a tool for studying the question, see \cite{rw3}: Under what conditions on $(M, \calf, g_0)$ the EGF metrics $g_t$
converge to one for which $\calf$ is umbilical, geodesic, or minimal\,?

Here, just for the readers' convenience, we provide the following definition from foliation theory, see \cite{cc}.

\begin{df}\rm
A family ${\mathcal F}=\{L_\alpha\}_{\alpha\in A}$ of connected
subsets of manifold $M^m$ is said to be an $n$-\textit{dimensional foliation}, if

1)\ $\bigcup_{\alpha\in A}L_\alpha=M^m$,

2)\ $\alpha\ne\beta\Rightarrow L_\alpha\bigcap L_\beta=\emptyset$,

3)\ for any point $p\in M$ there exists a $C^r-$chart (local coordinate system)
$\varphi_q: U_q\to\RR^m$ such that $q\in U_q$, $\varphi_q(q)=0$, and if
$\ U_q\bigcap L_\alpha\ne\emptyset\ $ the connected components of the sets $\varphi_q(U_q\bigcap L_\alpha)$ are given by
equations
$x_{n+1} = c_{n+1},\ \ldots,\ x_m=c_m,$
where $c_j$'s are constants.
The sets  $L_\alpha$ are immersed submanifolds of $M$ called {\it leaves} of $\calf$.
The family of all the vectors tangent to the leaves is an integrable subbundle of $TM$ denoted by $T\calf$.
If $M$ carries a Riemannian structure, $T\calf^\perp$ denotes the subbundle of all the vectors orthogonal to the leaves.
A foliation $\calf$ is said to be {\it orientable} (resp., {\it transversely orientable})
if the bundle $T\calf$ (resp. $T\calf^\perp$) is orientable.
\end{df}

\underline{Throughout the paper},
 $(M^{n+1}, g)$ is a Riemannian manifold with a codimension one transversely oriented foliation $\calf$,
$\nabla$ the Levi-Civita connection of~$g$, $N$ the positively oriented unit normal to $\calf$,
$A:X\in T\calf\mapsto-\nabla_X N$ the Weingarten opera\-tor of the leaves,
which we extend to a $(1,1)$-tensor field on $TM$ by $A(N)=0$.

Denote $\,\hat{}\,$ the $\calf$-component of a vector.
The definition $\hat S(X,Y):=S(\hat X, \hat Y)$ of the $\calf$-{\it truncated} $(0,2)$-{\it tensor} $S$ will be helpful throughout the paper.

Let $b:T\calf\times T\calf\to\RR$ be the second fundamental form of (the leaves of)
$\calf$ with respect to  $N$, and $\hat{b}$ its extension to
the $\calf$-truncated symmetric $(0, 2)$-tensor field on $M$.
Notice that $\hat{b}(N, \cdot)=0$ and $\hat{b}(X,Y)=g(A(X),Y)$.

The \textit{power sums} of the principal curvatures $k_1, \ldots, k_n$
of the leaves of $\calf$ (the eigenvalues of $A$) are given~by
\begin{equation*}
 \tau _j = k_1^j + \ldots + k_n^j = \tr(A^j)\quad (j\ge0).
\end{equation*}
They can be expressed using \textit{elementary symmetric functions}
$\sigma_0, \dots, \sigma_n$
$$ \sigma _j = \sum\nolimits_{i_1 <\ldots <i_j} k_{i_1}\cdot\ldots\cdot
k_{i_j}\quad (0\le j\le n),$$
called \textit{mean curvatures} in the literature (see Section~\ref{sec:gencomp}).

\begin{rem}\label{R-newton}\rm Certainly, the functions $\tau_{n+i}\ (i>0)$, are not independent:
they can be expressed as polynomials of $\vec\tau=(\tau _1,\ldots ,\tau _n)$,
using the Newton formulae, which in matrix form look like
\begin{equation*}
 T_n\left(\begin{smallmatrix}
 \sigma_1 \\ \sigma_2 \\ \dots\\ \sigma_n \\
 \end{smallmatrix}\right) =
 \left(\begin{smallmatrix}
 \tau_1 \\ \tau_2 \\ \dots\\ \tau_n \\
\end{smallmatrix}\right),
 \ \ \mbox{where}\ \
 T_n=\left(\begin{smallmatrix}
 1 &  &  & & & \\ \tau_1 & -2 &  &  &   \\
 \dots & \dots & \dots &  \dots &   \\
 \tau_{n-1} & \ -\tau_{n-2}  &\ldots &
 \ (-1)^{n}\tau_1 & \ (-1)^{n+1}n \\
 \end{smallmatrix}\right).
 \end{equation*}
Hence $\sigma$'s can be expressed by $\tau$-s using $T_n^{-1}$.
Moreover, we have the identities
\begin{eqnarray*}
 \tau_i=\det\left(\begin{smallmatrix}
 \sigma_1 & 1 & 0 & \ldots &  \ldots\\
 2\sigma_2& \sigma_1 & 1 & 0 & \ldots  \\
 \cdots & \cdots & \cdots & \cdots &  \\
 \ i\,\sigma_i \ \ & \ \sigma_{i-1} \ \ & \ \ldots \ \ & \ \sigma_2 \ \ & \sigma_1 \\
 \end{smallmatrix}\right),\quad i=1,2,\ldots
\end{eqnarray*}
Both formulae can be used to express $\tau_{n+i}$'s as polynomials of $\tau_{1},\ldots,\tau_{n}$.
\end{rem}

We study two types of evolution of Riemannian structures,
depending of functions $f_j\ (0\le j<n)$ in use (at least one of them is not identically zero):

(a) $f_j\in C^{2}(M\times\RR)$,
\ and

(b) $f_j= \tilde{f}_j(\vec\tau , \cdot)$, where $\tilde{f}_j\in C^{2}(\RR^{n+1})$.

Sometimes we will assume that $f_j=\tilde{f}_j(\vec\tau)$ with $\tilde{f}_j\in C^{2}(\RR^{n})$.

\begin{df}\rm
Given functions $f_j$ of type either (a) or (b),
a family $g_t,\ t\in[0,\eps)$, of Riemannian structures on $(M,\calf)$
will be called an \textit{Extrinsic Geometric Flow} (EGF), whenever
\begin{equation}\label{eq1}
 \dt g_t = h(b_t), \quad\text{where} \
 h(b_t)=h_t=\sum\nolimits_{j=0}^{n-1} f_j\,\hat{b}_j^{\,t}.
\end{equation}
Here, $\hat b_j^{\,t}$ are symmetric $(0,2)$-tensor fields on $M$ $g_t$-dual to $(A_t)^j$.
\end{df}

Hence, the EGF is an evolution equation which deforms Riemannian metrics by evolving them along $\calf$ in the direction of the tensor $h(b_t)$.
The~EGF preserves $N$ unit and perpendicular to $\calf$,
the $\calf$-component of the vector does not depend on $t$.
 One can interpret EGF as integral curves of a vector field $g\mapsto h(b)$, $h(b)$
being the right hand side of (\ref{eq1}), on the space $\mathcal{M}^k (M) = C^k(M, S_2^+ (M))$ of Riemannian
C$^k$-structures on $M$. Here, certainly, $S_2^+(M)$ is the bundle of positive definite
symmetric $(0, 2)$-tensors on~$M$.
This vector field may depend on time or not.

The choice of the right hand side in (\ref{eq1}) for $h(b_t)$ seems to be
natural: the powers $\hat{b}_j$ are the only (0,2)-tensors which can be
obtained algebraically from the second fundamental form $b$,
while $\tau_1, \ldots, \tau_n$ (or, equivalently, $\sigma_1, \ldots, \sigma_n$)
generate all the scalar invariants of extrinsic geometry. The powers
$\hat{b}_j$  with $j > 1$ are meaningful: for example,
the EGF with $h$ produced by

-- the extrinsic Ricci curvature tensor $\Ric^{\rm\,ex}(b)$, see (\ref{E-ricA2}),
depends on $\hat{b}_1,\hat{b}_2$;

-- the extrinsic Newton transformation $T_i(b)=\sigma_i\,\hat g-\sigma_{i-1}\,\hat b_1+\ldots+(-1)^i\hat b_i$

$(0<i<n)$, see Section~\ref{sec:appl}, depends on all $\hat{b}_j\ (1\le j\le i)$.

\vskip1mm
In this paper we investigate existence/uniqueness and properties of metrics $g_t$ satisfying (\ref{eq1}),
discuss in more details some particular cases and study some examples.
In Section~\ref{sec:main}, we collect main results.
Section~\ref{sec:pril} contains auxiliary results.
In Section~\ref{sec:gen-case}, we prove main results.
Section~\ref{sec:appl} is devoted to particular cases and examples.
Studying EGF will be continued in a series of papers remaining in preparation (see \cite{rw1}, \cite{rw3}).

\section{Main results}
\label{sec:main}

We shall omit index $t$ for $t$-dependent tensors $A$, $b,\,\hat b_j$ and functions $\tau_i, \sigma_i$.

The following theorem concerns the EGF of type (a)
and is essential in the proof of Theorem~\ref{T-localsol3} (for the EGF of type (b)).

\begin{thm}\label{L-main1}
Let $(M, g_0)$ be a compact Riemannian manifold with a codimen\-sion-1 foliation $\calf$.
Given functions $\tilde f_j\in C^{2}(M\times\RR)$
there is a unique smooth solution $g_t$ to (\ref{eq1}) of type (a) defined on some positive time interval $[0,\eps)$.
\end{thm}

In particular, there exists a unique smooth solution $g_t,\ t\in[0,\eps)$ to
\[
 \dt g_t=\sum\nolimits_{j=0}^{n-1} a_j(t)\,\hat b_j,\quad a_j\in C^{2}(\RR).
\]
Although (\ref{eq1}) of type (b) consists of first order non-linear PDEs,
the corresponding power sums $\tau_i\ (i>0)$ satisfy the infinite
quasilinear system
\begin{equation}\label{E-tauAk}
 \dt\tau_i+\frac{i}2\big\{\hskip-1pt\tau_{i-1}N(f_0)
 +\hskip-2pt\sum\nolimits_{j=1}^{n-1}\hskip-3pt\big[\frac{j f_j}{i{+}j{-}1}
 N(\tau_{i+j-1}) +\tau_{i+j-1}N(f_j)\big]\hskip-1pt\big\} = 0,
\end{equation}
where $N(f_j)=\sum\nolimits_{s=1}^n
f_{j,\tau_s}\,N(\tau_s)$. By Proposition~\ref{L-BBnm}
(Section~\ref{sec:gencomp}), the matrix of the $n$-truncated system
(\ref{E-tauAk}) (where $\tau_{n+i}$'s are replaced by suitable polynomials of
$\tau_1,\ldots,\tau_n$) has the form $\tilde A + \tilde B$, where $\tilde A = (\tilde{A}_{ij})$,
\begin{equation}\label{E-ABnm}
 \tilde A_{ij}=\frac i2\sum\nolimits_{m=0}^{n-1} \tau_{i+m-1} f_{m,\tau_j},\ \
 \tilde B=\sum\nolimits_{m=1}^{n-1}\frac m2\,f_m\cdot(B_{n,1})^{m-1},
\end{equation}
and $B_{n,1}$ is the \textit{generalized companion matrix} (see Section~\ref{sec:gencomp}) to
the characteristic polynomial of $A_N$.
 Recall that an $n$-by-$n$ matrix is \textit{hyperbolic} (see Section~\ref{sec:quasi})
if its right eigenvectors are real and span $\RR^n$.

\begin{ex}\rm
(a) If $\dt g_t=f(\vec\tau,t)\,\hat g_t$ (i.e., $f_j=\delta _{j 0} f$),
then, by ~(\ref{E-ABnm}), $\tilde B=0$ and $\tilde A=(\tilde a_{ij})$, where
$\tilde a_{ij}=(i/2)\,\tau_{i-1} f_{,\tau_j}(\vec\tau,0)$.
System (\ref{E-tauAk}) reduces to
\begin{equation}\label{E-PDE-1b-0a}
 \dt\tau_i + \frac i2\,\tau_{i-1}\sum\nolimits_{j=1}^n f_{,\tau_j}(\vec\tau,t)\,N(\tau_j)=0,
 \quad i=1,2,\ldots
\end{equation}
The matrix $\tilde A$ of $n$-truncated system (\ref{E-PDE-1b-0a}) is hyperbolic if
for any $p\in M$

\vskip3mm
 {\rm either}
 $2\tr\tilde A\,{=}\!\sum\limits_{i} i\,\tau_{i-1} f_{,\tau_i}\,{\ne}\,0$
 {\rm or}
 $\tilde A\equiv 0$ on the $N$-curve through $p$.
 \hfill $(H_1)$

\vskip2mm\noindent
(b) If $\dt g_t=f(\vec\tau,t)\,\hat b_1$ (i.e., $f_j=\delta_{j 1} f$)
then, again by~(\ref{E-ABnm}), $\tilde B=\frac12\,f(\vec\tau,0)\id$ and $\tilde A=(\tilde a_{ij})$,
where $\tilde a_{ij}=(i/2)\,\tau_{i}\,f_{,\tau_j}(\vec\tau,0)$.
System (\ref{E-tauAk}) reduces to
\begin{equation}\label{eq322}
 \dt\tau_i +\frac12 f(\vec\tau,t) N(\tau_i)
 +\frac i2\,\tau_{i}\sum\nolimits_{j=1}^n f_{,\tau_j}(\vec\tau ,t)\,N(\tau_j)=0,\quad i=1,2,\ldots
\end{equation}
The matrix $\tilde A+\tilde B$ of $n$-truncated system (\ref{eq322}) is hyperbolic if
for any $p\in M$

\vskip3mm
 {\rm either}
 $2\tr\tilde A=\sum\limits_{i} i\,\tau_{i} f_{,\tau_i}\ne 0$
 \ {\rm or}
 $\tilde A\equiv0$ on the $N$-curve through $p$.
 \hfill $(H_2)$
\end{ex}

 The central result of our paper is the following.

\begin{thm}[\textbf{Short time existence}]\label{T-localsol3}
Let $(M, g_0)$ be a compact Riemannian manifold with a codimension-one foliation $\calf$ and a unit normal $N$.
If~the matrices $\tilde A+\tilde B$ and $\tilde B$ of (\ref{E-ABnm}) are hyperbolic for all $p\in M$ and $t=0$, then the EGF (\ref{eq1}) of type (b) has a unique smooth solution $g_t$ defined on some positive time interval $[0,\eps)$.
\end{thm}

The proof of Theorem~\ref{L-main1} follows standard methods of the theory of first order PDEs with one space variable. The proof of Theorem~\ref{T-localsol3} (Section~\ref{sec:Texist}) consists of several steps.

1) The power sums $\tau_i$ are recovered on $M$ (as a unique solution to a quasilinear hyperbolic system of PDEs) for some positive time interval $[0,\eps)$, see Lemmas~\ref{L-Cn2}, \ref{L-symm-k1} and \ref{P-tauAk}, and Proposition~\ref{L-BBnm}.

2) Given $(\tau_i)$ (of Step 1), the metric $g_t$ is recovered on $M$
(as a unique solution to certain quasilinear system of PDEs),
 see Theorem~\ref{L-main1} and Lemma~\ref{L-adapted-coords}.

3) The $\tau_i$-s of the $g_t$-principal curvatures of $\calf$ (of Step 2) are shown
to coincide with $\tau_i$ (of Step 1), see Theorem~\ref{L-main1} and Lemmas~\ref{P-Ak} and \ref{L-Tlocalsol2}.

\vskip1mm
Let us remark that the solution in Theorem~\ref{T-localsol3} is unique if only $\tilde A+\tilde B$ is hyperbolic.
For $f_j=0\ (j\ge2)$, Theorem~\ref{T-localsol3} holds under weaker condition
that only the matrix $\tilde A$ is hyperbolic for all $p\in M$ and $t=0$.

Denote by $\mathcal{L}_Z$ the Lie derivative along a vector field $Z$.

\begin{cor}\label{C-A1}
Let $(M, g_0)$ be a compact Riemannian manifold with a codimen\-sion-1
foliation  $\calf$ and a unit normal $N$.
If $f\in C^{2}(\RR^{n+1})$ and the condition $(H_2)$ is satisfied at $t=0$ and any $p\in M$,
then there exists a unique smooth solution $g_t,\ t\in[0,\eps)$ to the EGF \begin{equation}\label{eq-g-f1}
 \dt g_t=f(\vec\tau, t)\,\hat b_1,\quad t\in[0,\eps)
\end{equation}
for some $\eps>0$. Furthermore, $g_t$ can be determined from the system
 $\mathcal{L}_{Z_t}\, g_t =0$
 with
 $Z_t=\dt+\frac12f(\vec\tau,t)\,N$,
where $\vec\tau$ are the unique smooth solution to (\ref{eq322}).
\end{cor}

\begin{cor}\label{C-A0}
Let $(M, g_0)$ be a compact Riemannian manifold with a codimen\-sion-1 foliation
$\calf$ and a unit normal $N$.
If $f\in C^{2}(\RR^{n+1})$ and the condition $(H_1)$ is satisfied
at $t=0$ and any $p\in M$,
then there exists a unique smooth solution $g_t$ to the EGF \begin{equation}\label{E-gC-A0}
 \dt g_t=f(\vec\tau,t)\,\hat g_{t},\quad t\in[0,\eps)
\end{equation}
for some $\eps>0$. Furthermore,
 $\hat g_t= \hat g_0\exp(\int_0^t f(\vec\tau,t)\,dt)$,
where the power sums $\vec{\tau}$ are the unique solution to (\ref{E-PDE-1b-0a}).
\end{cor}

In the particular case $f_1 = C = const$, $f_j = 0\ (j\ne 1)$, the system (\ref{E-tauAk}) (see also (\ref{eq322}) for $f=C$) is reduced to the linear PDE
\begin{equation*}
 \dt\tau_i +(C/2) N(\tau_i)=0.
\end{equation*}
The above equation can be interpreted on $M\times\RR$ by saying that
$\tau_i$ is constant along the orbits of the vector field
$X = \dt +(C/2)N$. If $(\psi_t)$ denotes the flow of $(C/2)N$ on $M$, then the flow
$(\phi_t)$ of $X$ is given by $\phi_t(p, s) = (\psi_t(p), t + s)$ for
$p\in M$ and $s\in\RR$, therefore $\phi_t$ maps the level surface
$M_s = M\times\{s\}$ onto $M_{t+s}$, in particular, $M_0$ onto $M_t$.
This implies

\begin{cor}\label{cor1-tau}
If $f_1 = const$ and $f_j = 0$ for all $j\ne 1$ (for the EGF), then for all $i$ and $t$, the following equality holds:
 $\tau_i^t = \tau_i^0\circ\phi_{-t}\,$.
In particular, if $\tau_i^0={\rm const}$ for some $i$, then $\tau_i={\rm const}$ for all $t$.
\end{cor}

Recall that a compact manifold $M$ equipped with a codimension-one
foliation $\calf$ admits a Riemannian structure $g$ for which all the
leaves are minimal ($\tau_1=0$ in the above terminology) if and only if
$\calf$ is {\it topologically taut}, that is every its leaf meets a
loop transverse to the foliation \cite{su}. The known proofs of existence of
such  metrics use Hahn-Banach Theorem and do not show how to construct
them. Above observations show how to produce a $1$-parameter family of
metrics with $\tau_{\,2j+1}=0$ (with fixed $j$) starting from one of them.

\vskip1mm
For multi-index
 ${\bm\alpha}=(\alpha_1,\ldots,\alpha_n)\in\ZZ^n_+$
define
 $\vec\tau^{\,\bm\alpha}:=\tau_1^{\alpha_1}\ldots\tau_n^{\,\alpha_n}.$
Set
 $J_{m,l}=\{{\bm\alpha}\in\ZZ^n_+: \sum\nolimits_j\alpha_j=m, \ \sum\nolimits_j j\,\alpha_j=l\}.$

Recall that a vector field on a manifold $M$ is \textit{complete},
if any its trajectory $c: t\mapsto c(t)$ can be extended to the whole range $\RR$ of parameter~$t$.
If $M$ carries a Riemannian structure $g$ and a vector field $X$ has bounded length,
then the completeness of $(M, g)$ is sufficient for the completeness of $X$.

\begin{prop}\label{T-f1-B}
Let $(M, g_0)$ be a Riemannian manifold with a codimen\-sion-1 foliation $\calf$ and a complete unit normal field $N$.
Given $c_{\bm\alpha}\in\RR\ ({\bm\alpha}\in J_{m,l})$ and $m,l\in\NN$, define functions
$f_t=\sum\nolimits_{\,{\bm\alpha}\in J_{m,l}}c_{\bm\alpha}\,(\vec\tau^{\,t})^{\bm\alpha}$
(where $\vec\tau^{\,t}$ is, as usual, the vector of power sums of principal curvatures of the leaves)
and set

\vskip1.5mm
 $T=\infty$ \ if $N(f_0)\ge0$ on $M$ \ and $T=-\frac{2}{l+1}/\inf_{\,M} N(f_0)$ \ otherwise.

\vskip1.5mm\noindent
Then the EGF $\dt g_t=f_t\,\hat b_1$, compare (\ref{eq-g-f1}), has a unique smooth solution on $M$ for $t\in[\,0, T)$ and does not possess one for $t\in[\,0, T]$.
\end{prop}

We apply EGF to \textit{totally umbilical} foliations
(that is, such that the Weingarten operator $A$ is proportional to the identity at any point,
among them \textit{totally geodesic} foliations appear when $A=0$) and to foliations on surfaces.
First, we shall show that EGF preserve total umbilicity of $\calf$.

\begin{prop}\label{P-08}
Let $(M, g_0)$ be a Riemannian manifold endowed with a codi\-mension-1 totally umbilical foliation  $\calf$. If
$g_t\ (0\le t<\eps)$ provide the EGF of type (b) on $(M,\calf)$,
then $\calf$ is totally umbilical for any $g_t$.
\end{prop}

\begin{prop}\label{C-tu1}
Let $(M, g_0)$ be a Riemannian manifold,
and $\calf$ a codimen\-sion-1 totally umbilical foliation on $M$ with the normal curvature
$\lambda_0$ and a complete unit normal field~$N$. Set

\vskip1.5mm
$T=\infty$ if $N(\psi(\lambda_0))\ge0$ on $M$, and $T=-2/\inf_{M}N(\psi(\lambda_0))$ otherwise,

\vskip1.5mm\noindent
where $\psi (\lambda) =\sum\limits_{j=0}^{n-1} f_j(n\lambda ,\ldots ,n\lambda ^n)\,\lambda ^j$.
Then the EGF with $h=\sum\limits_{j=0}^{n-1} f_j(\vec\tau)\,\hat{b}_j$
has a unique smooth solution $g_t$ on $M$
for $t\in[0, T)$, and does not possess one for $t\ge T$.
Moreover, $\calf$ is $g_t$-totally umbilical and
$\hat g_t= \hat g_0\exp(\int_0^t\psi(\lambda_t)\,dt),$
where $\lambda _t$ is a unique smooth solution to the PDE
\begin{equation}\label{E-tauAk-L}
 \dt\lambda_t+\frac12\,N(\psi'(\lambda_t))=0.
\end{equation}
\end{prop}

\section{Preliminaries}
\label{sec:pril}

\subsection{Hyperbolic quasi-linear PDEs}
\label{sec:quasi}

Here, we recall known results on quasi-linear PDEs.

Let $A=(a_{ij}(x,t,\bu))$ be an $n\times n$ matrix,
$\vec b=(b_{i}(x,t,\bu))$ -- an $n$-vector.
A~first order \textit{quasilinear} system of PDEs, $n$ equations in $n$ unknown functions $\bu = (u_1,\ldots ,u_n)$ and two variables $x,t\in\RR$, has the form
\begin{equation}\label{E-PDE-1}
 \dt\bu + A(x,t,\bu)\,\dx\bu = \vec b(x,t,\bu).
\end{equation}
When the coefficient matrix $A$ and the vector $\vec b$ are functions of $x$ and $t$ only, the system is just \textit{linear};
if $\vec b$ only depends also on $\bu$, the system is said to be \textit{semilinear}.
The \textit{initial value problem} for (\ref{E-PDE-1}) with initial surface
$\Pi =\{t = 0\}$ and given smooth data $\bu_0,\ A$ and $\vec b$ consists in finding smooth $\bu$ such that (\ref{E-PDE-1}) and $\bu(x,0)=\bu_0(x)$ are satisfied.

\begin{df}\rm
 The system (\ref{E-PDE-1}) is \textit{hyperbolic} in the $t$-direction at $(x,t,\bu)$
(in an appropriate domain of the arguments of $A$ and $\vec b$) if the right eigenvectors of $A$ are real and span $\RR^n$.
For a solution $\bu(x,t)$ to (\ref{E-PDE-1}), the corresponding eigenvalues $\lambda_i(x,t,\bu)$ are called the \textit{characteristic speeds}.
 The system is \textit{strictly hyperbolic} if $\lambda_i(x,t,\bu)$ are distinct.
 For the hyperbolic system (\ref{E-PDE-1}), the vector field
 $\partial_t+\lambda_i(x,t,\bu)\partial_x$ is called the $i$-\textit{characteristic field}, and its integral curves are called $i$-\textit{characteristics}.
\end{df}

\noindent
Remark that the \textit{hyperbolicity of} $A$ is equivalent to any of the properties:

``$A$~has real eigenvalues $\lambda_1\le\ldots\le\lambda_n$ and simple elementary divisors"

 (i.e., $A$ has no Jordan cells of order greater than one), and 

''$A$ is diagonal in some affine basis".

\noindent
 Hence, the hyperbolic matrix $A$ can be represented as $A = R D R^{-1}$, where $R$ is a nonsingular $n \times n$ matrix and $D$ is a diagonal matrix. The columns $r_i$ of $R$ are the right eigenvectors of $A$, whereas the rows of $R^{-1}$ are left eigenvectors of $A$.
A hyperbolic system reduces to the ODEs for its characteristic fields.
Indeed, multiplying (\ref{E-PDE-1}) by $\vec r_i^{\ T}$ and using $\frac{d}{dt}\bu=\dt\bu + \lambda_i\,\dx\bu$, we get the~ODE
\[
 r_i^{\ T}\,{d}\bu/{dt}=r_i^{\ T}\,b\ \
 \mbox{\rm along the characteristic } \ {d}x/{dt}=\lambda_i(x,t,\bu).
\]

\textbf{Theorem A} \cite{hw}
\textit{\,Let the quasi-linear system (\ref{E-PDE-1}) be such that}

\hskip-2mm
1) \textit{it is hyperbolic in the $t$-direction in $\Omega=\{|x|\le a,\ 0\le t\le s,\ \|\vec u\|_{\infty}\le r\}$}

   \textit{for some $s,r>0$},

\hskip-2mm
2) \textit{the matrix $A$ and the vector $\vec b$ are $C^{1}$-regular in $\Omega$}.

\noindent
 \textit{If the initial condition
\begin{equation}\label{E-pdes-2}
  \bu(0,x)= \bu_0(x),\quad x\in [-a, a]
\end{equation}
has $C^{1}$-regular $\bu_0$  in $[-a,a]$ and $\|\bu_0\|_{\infty}<r$,
then (\ref{E-PDE-1})\,--\,(\ref{E-pdes-2}) admit a unique $C^{1}$-regular solution $\vec u(x,t)$
in $\bar\Omega=\{(x,t):\ |x| + K t\le a,\ 0\le t\le\eps\}$,
with $K=\max\limits\{|\lambda_i(x,t,\vec u)|:\ (x,t,\vec u)\in\Omega,\ 1\le i\le n\}$.}

\begin{ex}\label{Ex-burgers}\rm
For any function $\psi\in C^{1}(\RR)$, we can multiply the equation
\begin{equation}\label{E-PDE-1-c1}
  \dt u + \psi(u)\,\dx u = 0,
\end{equation}
 across by $\psi\,'(u)$, and obtain
 $\dt \psi + \psi\cdot\dx\psi = 0$ (the Burgers' equation).
Thus, the behaviour of the solutions to (\ref{E-PDE-1-c1}) (for $t$ before the first singular value)
is not expected to be much different from that for Burgers' equation.
\end{ex}

\subsection{The generalized companion matrix}
\label{sec:gencomp}

Let $P_n=\lambda^n-p_1\lambda^{n-1}-\ldots-p_{n-1}\lambda-p_n$ be a polynomial over $\RR$
and $\lambda_1\le\lambda_2\le\ldots\le\lambda_n$ be the roots of $P_n$ for $n>1$.
Hence, $p_i=(-1)^{i-1}\sigma_i$, where $\sigma_i$ are elementary symmetric functions of the roots $\lambda_i$.

\begin{df}\rm
The \textit{generalized companion matrices} of $P_n$ are defined as
$$
 \hat C_{g}=\hskip-.3mm\left(\hskip-.5mm\begin{smallmatrix}
                     0 & \frac{c_{n{-}1}}{c_n} & 0 & \cdots & 0\\
                     0 & 0 & \frac{c_{n{-}2}}{c_{n{-}1}} & \cdots  & 0\\
                     \cdots & \cdots & \cdots & \cdots & \cdots\\
                     0 & 0 & \cdots & 0 & \frac{c_{1}}{c_{2}} \\
           c_{n} p_n& c_{n{-}1}p_{n{-}1}&\ldots & c_{2} p_2& c_{1} p_1\\
     \end{smallmatrix}\hskip-.5mm\right)\hskip-.5mm
       \quad\mbox{or}\quad
 \check{C}_{g}=\hskip-.3mm\left(\hskip-.5mm\begin{smallmatrix}
           c_{1} p_1& c_{2} p_{2}&\ldots& c_{n{-}1} p_{n{-}1}& c_{n} p_n\\
                     \frac{c_{1}}{c_2} & 0 & \cdots & 0 & 0\\
                     0 &\frac{c_{2}}{c_{3}} & 0 & \cdots  & 0\\
                     \cdots & \cdots & \cdots & \cdots & \cdots\\
                     0 & \cdots & 0 & \frac{c_{n-1}}{c_{n}} & 0 \\
     \end{smallmatrix}\hskip-.5mm\right),
$$
where $c_1=1$ and $c_i\ne0\ (i>1)$ are arbitrary numbers.
\end{df}

Notice that $\hat C_g$ acts on $\RR^n$ as $\bx\to\hat C_g\bx$, where $\bx{=}(x_1,\ldots,x_n)$.
Inverting the order of indices, i.e., taking $(x_n,\ldots,x_1)$, one may describe this action by~$\check C_{g}$.
 If  all $c_i$'s are equal to $1$, the matrix $\hat C_{g}$
reduces to the standard \textit{companion matrix} $C_{g}$ of $P_n$.
The explicit formulae (polynomials) for entries in powers of $C_{g}$
and some applications to the theory of the symmetric functions are given in~\cite{cl}.
The following matrix $B_{n,1}$ (the \textit{generalized companion matrix} with $c_i=\frac{n}{n-i+1}$) plays the key role in the paper:
\begin{equation}\label{E-Bn1}
 B_{n,1}=\left(\begin{smallmatrix}
           0 & \frac12 & 0 & \cdots & 0\\
           0 & 0 & \frac23 & \cdots  & 0\\
           \cdots & \cdots & \cdots & \cdots & \cdots\\
           0 & 0 & 0 & \cdots & \frac{n-1}{n} \\
           (-1)^{n-1}\frac{n}{1}\,\sigma_n\quad& (-1)^{n-2}\frac{n}{2}\,\sigma_{n-1}\quad& \ldots\quad& \ -\frac{n}{n-1}\,\sigma_2\ \ & \sigma_1 \\
     \end{smallmatrix}\right).
\end{equation}

\begin{lem}\label{L-Cn2}
The generalized companion matrices have the properties

\hskip-12pt
a) the characteristic polynomial of $\hat C_g$ (or $\check{C}_g$) is
$P_n$.

\hskip-12pt
b) $v_j=(1, \frac{c_{n}}{c_{n-1}}\lambda_j,\frac{c_{n}}{c_{n-2}}\lambda_j^2,\ldots,c_{n}\lambda_j^{n-1})$
is the eigenvector of $\hat C_g$ for the eigenvalue $\lambda_j$,

 resp., $w_j=(c_{n}\lambda_j^{n-1},\ldots,\frac{c_{n}}{c_{n-2}}\lambda_j^2,\frac{c_{n}}{c_{n-1}}\lambda_j, 1)$
 is the eigenvector of $\check{C}_g$ for $\lambda_j$.

\hskip-12pt
c) $\hat C_g V = V D$, where $V\!=\!\{\frac{c_{n}}{c_{n-i+1}}\lambda_j^{i-1}\}_{1\le i,j\le n}$ is a Vandermonde type matrix,

  and $D={\rm diag}(\lambda_1,\ldots,\lambda_{n})$ a diagonal matrix.
  (If all $\lambda_i$'s are distinct,

  then obviously $V^{-1}\hat C_g V =  D$).
 \end{lem}

\noindent\textbf{Proof}.
a) We will show by induction for $n$ that
 $\det|\lambda\id_n-\hat C_{g}|=P_n$,
hence the eigenvalues $\lambda_{1}\le\ldots\le\lambda_{n}$ of $\hat C_{g}$ coincide with the roots of $P_n$.

Expanding by co-factors down the first column, we obtain
$$
 \det|\lambda\id_n-\hat C_{g}|=
 \lambda\,P_{n-1}-(-1)^{n-1}c_n\,p_n\prod\nolimits_{i=1}^{n-1}\big(-{c_i}/{c_{i+1}}\big),
$$
where
 $P_{n-1}=\lambda^{n-1}-p_1\lambda^{n-2}-\ldots-p_{n-2}\lambda-p_{n-1}$ (by the induction assumption) is a certain polynomial of degree $n-1$. Notice that $c_n\prod_{i=1}^{n-1}\frac{c_i}{c_{i+1}}{=}\,1$.
This completes the proof of the claim.

b) The direct computation shows us that
$$(\lambda_j \id_n-\hat C_{g})\,\vec v_j=0,\qquad
   (\lambda_j \id_n-\check{C}_{g})\,\vec w_j=0.
$$

c) Hence $\hat C_g V = V D$, where $D={\rm diag}(\lambda_1,\ldots,\lambda_{n})$
is a diagonal matrix.
If  $\{\lambda_j\}$ are pairwise distinct, then $\det V\ne0$ and obviously
$V^{-1}\hat C_{g}V=D$.$\hfill\square$

\vskip1.5mm
Consider the infinite system of linear PDEs with functions $f_j\in C^{2}(\RR^2)$
\begin{equation}\label{E-tauBnm}
  \dt\tau_i+\frac{i}2\sum\nolimits_{j=1}^{n-1}\frac{j}{i{+}j{-}1}\,f_j(t,x)\,\dx\tau_{i+j-1}=0,\qquad i=1,2,\ldots
\end{equation}
where $\tau _i\ (i\in\mathbb{N})$ are the power sums of smooth functions $\lambda _i(t,x)\ (i\le n)$.
Let $\sigma _j\ (j\le n)$ be the elementary symmetric functions of $\{\lambda _i(t,x)\}$.

\begin{prop}\label{L-BBnm}
The matrix of $n$-truncated system (\ref{E-tauBnm})
(where $\tau_{n+i}$'s are eliminated using suitable polynomials of $\tau _1,\ldots ,\tau _n$,
as described in Remark~\ref{R-newton}) is $\tilde B=\sum\nolimits_{m=1}^{n-1} f_m B_{n,m-1}$,
where
\begin{equation}\label{E-BBB}
 B_{n,m}=\frac{m+1}2\,(B_{n,1})^{m}
\end{equation}
and $B_{n,1}$ is the generalized companion matrix (\ref{E-Bn1}).
The eigenvalues of $\tilde B$ are
$$\tilde{\lambda}_j=\frac12\sum\nolimits_{m=1}^{n-1} m\,f_m\,\lambda_j^{m-1}.$$
while the corresponding eigenvectors are given by
$v_i=(1, 2\,\lambda_i, 3\,\lambda_i^2,\ldots,  n\,\lambda_i^{n-1})$.
\end{prop}

\noindent\textbf{Proof}. We need to prove (\ref{E-BBB}) only.
Replacing $\dx\tau_{n+j}$ in (\ref{E-tauBnm}) by linear combinations
of $\dx\tau_{i}\ (i\le n)$ due to (\ref{E-dtaunm})
in what follows, we find the $(i,j)$ entry of $B_{n,m}$
\begin{equation}\label{E-bij}
 b_{ij}^{(n,m)}=
 \Big\{\begin{array}{cc}
 \frac{i\,(m+1)}{2(i+m)}\,\delta^j_{i+m} & \mbox{ if } \ i+m\le n,\\
 (-1)^{n-j}\,\frac{i\,(m+1)}{2 j}\,\beta_{n,\,i+m-n,\,j} & \mbox{ if } \ i+m > n.\\
 \end{array}
\end{equation}
In particular, $B_{n,0}=(1/2) \id_n$.
The equality (\ref{E-BBB}) follows directly from
\[
 B_{n,m}=\frac{m+1}{m}B_{n,m-1}B_{n,1}.
\]
Notice that the last formulae are true for $m=1$.
We shall show that all the $(i,j)$-entries of the matrices $\frac{m+1}{m}B_{n,m-1}B_{n,1}$
and $B_{n,m}$ coincide.

First, let  $i+m-1\le n$. Then,
\begin{eqnarray*}
 \frac{m+1}{m}\sum\nolimits_s b_{is}^{(n,m-1)} b_{sj}^{(n,1)}
 \overset{(\ref{E-bij})}=\frac{m+1}{m}\sum\nolimits_s\frac{im}{2(i+m-1)}\,
 \delta^s_{i+m-1}\frac{j-1}{j}\,\delta^s_{j-1}\\
 \overset{j=i+m}=\frac{m+1}{2}\frac{i}{i+m-1}\,\delta^{j-1}_{i+m-1}\frac{i+m-1}{i+m}
 =\frac{i(m+1)}{2(i+m)}\,\delta^{j}_{i+m}\overset{(\ref{E-bij})}=b_{ij}^{(n,m)}.
\end{eqnarray*}

Now, let $\tilde m = i+m-1-n > 0$. Then
\begin{eqnarray*}
 \frac{m+1}{m}\sum\nolimits_{s=1}^{\,n} b_{is}^{(n,m-1)} b_{sj}^{(n,1)}
 =\frac{m+1}{m}[b_{i,j-1}^{(n,m-1)} b_{j-1,j}^{(n,1)}
 +b_{i n}^{(n,m-1)}b_{nj}^{(n,1)}]\\
 \overset{(\ref{E-bij})}=\frac{m{+}1}{m}\big[(-1)^{n-j+1}\frac{im}{2(j-1)}\beta_{n,\tilde m,j-1}\frac{j{-}1}{j}
 +\frac{im}{2n}\beta_{n,\tilde m,n}(-1)^{n-j}\frac{n}{j}\sigma_{n-j+1}\big]\\
\overset{(\ref{E-bbeta})}
 =\frac{i(m+1)}{2 j}(-1)^{n-j}
 (\beta_{n+1,\tilde m,n+1}\sigma_{n-j+1}-\beta_{n+1,\tilde m,j})\\
 =\frac{i(m+1)}{2 j}(-1)^{n-j}\beta_{n,\tilde m+1,j}
 \overset{(\ref{E-bij})}=b_{ij}^{(n,m)}
\end{eqnarray*}
 By (\ref{E-BBB}) and (\ref{E-Bn1}), if $\lambda _j$'s are pairwise different,
 then the matrix $B_{n,m}$ $(m>0)$ is hyperbolic.$\hfill\square$

\begin{lem}\label{L-symm-k1}
The coefficients $\beta_{n,m,i}$ of the decomposition
\begin{equation}\label{E-dtaunm}
 \frac{1}{n+m}\,\dx\tau_{n+m}=
 \sum\nolimits_{i=1}^n (-1)^{n-i}\,\frac{1}{i}\,\beta_{n,m,i}\,\dx\tau_i,
 \qquad m>0
\end{equation}
satisfy the following recurrence relations:
\begin{eqnarray}\label{Eq-beta}
 \beta_{n,1,i}\eq\sigma_{n-i+1},\quad
 \beta_{n,m,i}=\beta_{n+1,m-1,n+1}\sigma_{n-i+1}-\beta_{n+1,m-1,i}\quad (m>1),\qquad\\
 \label{E-bbeta}
 \beta_{n,m,i}\eq\beta_{n+j,\,m,\,i+j}
 \quad (1\le i\le n,\ m,j>0).
\end{eqnarray}
\end{lem}

\begin{rem}\rm
 In view of (\ref{E-bbeta}), relation (\ref{Eq-beta})$_2$  reduces to
\begin{equation*}
 \beta_{n,m,i}=\beta_{n,m-1,n}\sigma_{n-i+1}-\beta_{n,m-1,i-1}\quad (m>1).
\end{equation*}
For small values of $m$, $m=1,2$, we have from (\ref{E-dtaunm})
\begin{eqnarray}\label{E-symm-k1}
 \frac{1}{n{+}1}\,\dx\tau_{n+1}\eq\sum\nolimits_{i=1}^n\frac{(-1)^{n-i}}{i}\,\sigma_{n-i+1}\dx\tau_i,\\
 \label{E-symm-k2}
 \frac{1}{n{+}2}\,\dx\tau_{n+2}\eq \sum\nolimits_{i=1}^n\frac{(-1)^{n-i}}{i}
 (\sigma_{1}\/\sigma_{n-i+1}-\sigma_{n-i+2})\,\dx\tau_i.\qquad
\end{eqnarray}
By Proposition~\ref{L-BBnm}, the last row of $B_{n,1}$ (resp., of $B_{n,2}$) consists of the coefficients at $\partial_x\tau_i$'s on the RHS of (\ref{E-symm-k1}),
(resp., of  (\ref{E-symm-k2})) and so on.
\end{rem}

\noindent\textbf{Proof of Lemma~\ref{L-symm-k1}}.
Let $m=1$. Equality $\tau_{i}=\sum_{j}{\lambda_j}^{i}$ yields
$\frac{1}{i}\,\dx\tau_{i}=\sum_{j}{\lambda_j}^{i-1}\dx \lambda_j$.
Similarly, we find
\begin{equation*}
 \sum\nolimits_{i=1}^n\frac{(-1)^{n-i}}{i}\,\sigma_{n-i+1}\dx\tau_i=
 \sum\nolimits_{j}\dx \lambda_j\sum\nolimits_{i=1}^n (-1)^{n-i}\,\sigma_{n-i+1}\lambda_j^{i-1}.
\end{equation*}
Define the polynomial $P_n(x)=\lambda^n{-}\sigma_1(x)\lambda^{n-1}{+}\ldots{+}(-1)^n\sigma_n(x)$.
Since $\lambda_j(x)$ are the roots of $P_n$, we obtain
\begin{eqnarray*}
 \frac{1}{n+1}\,\dx\tau_{n+1}&-&\sum\nolimits_{i=1}^n\frac{(-1)^{n-i}}{i}\,\sigma_{n-i+1}\dx\tau_i\\
 \eq\sum\nolimits_{j}\big(\lambda_j^n-\sigma_1\lambda_j^{n-1}+\ldots +(-1)^n\sigma_n\big)\dx\lambda_j=0
\end{eqnarray*}
that proves (\ref{E-symm-k1}). Hence, $\beta_{n,1,i}=\sigma_{n-i+1}$.

In order to prove the recurrence relation in (\ref{Eq-beta}),
assume temporarily that $\lambda_{n+1}=\eps$, $\tilde n = n+1$ and $\tilde m = m-1$.
Hence, $\frac{1}{n{+}m}\,\dx\tau_{n+m}=\frac{1}{\tilde n{+}\tilde m}\,\dx\tau_{\tilde n+\tilde m|\,\eps=0}$.
Then we put $\eps=0$ and replace $\dx\tau_{n+1}(x)$ due to (\ref{E-symm-k1})
\begin{eqnarray*}
 && \hskip-10mm
 \frac{1}{\tilde n{+}\tilde m}\,\dx\tau_{\tilde n+\tilde m}=
 \sum\nolimits_{i=1}^{\tilde n}\frac{(-1)^{\tilde n-i}}{i}\,\beta_{\tilde n,\tilde m,i}\,\dx\tau_i \\
 \eq\beta_{n+1,m-1,n+1}\frac{1}{n+1}\,\dx\tau_{n+1}
 +\sum\nolimits_{i=1}^{n}\frac{(-1)^{n-i+1}}{i}\beta_{n+1,m-1,i}\,\dx\tau_i\\
 \eq\hskip-6mm\stackrel{\eps=\,0}{\phantom{=}} \sum\nolimits_{i=1}^n\frac{(-1)^{n-i}}{i}
 \big(\beta_{n+1,m-1,n+1}\,\sigma_{n-i+1}-\beta_{n+1,m-1,i}\big)\dx\tau_i
\end{eqnarray*}
that completes the proof of (\ref{Eq-beta}).
For $m=2$, we deduce from (\ref{Eq-beta})
$$\beta_{n,2,i}=\beta_{n+1,1,n+1}\,\sigma_{n-i+1}-\beta_{n+1,1,i}=\sigma_{1}\sigma_{n-i+1}-\sigma_{n-i+2}$$
that proves (\ref{E-symm-k2}).

Finally, we prove (\ref{E-bbeta}) by induction on $m$. For $m=1$ we have
$$\beta_{n+j,\,1,\,i+j}\overset{(\ref{Eq-beta})_1}=\sigma_{(n+j)-(i+j)+1}=\sigma_{n-i+1}
\overset{(\ref{Eq-beta})_1}=\beta_{n,1,i}.$$
Assuming (\ref{E-bbeta}) for $m-1$ we deduce it for $m$ using (\ref{Eq-beta})$_2$:
\begin{eqnarray*}
\beta_{n+j,\,m,\,i+j}&\overset{(\ref{Eq-beta})_2}=&
 \beta_{(n+j)+1,\,m-1,\,(n+j)+1}\sigma_{(n+j)-(i+j)+1}-\beta_{(n+j)+1,\,m-1,\,i+j}\\
 \eq \beta_{n+1,\,m-1,\,n+1}\sigma_{n-i+1}-\beta_{n+1,\,m-1,\,i}
 \overset{(\ref{Eq-beta})_2}=\beta_{n,m,i}
\end{eqnarray*}
that completes the proof of (\ref{E-bbeta}).$\hfill\square$

\begin{ex}\label{Ex-ex-ric2}\rm
For $f_j=\delta_{j1}$, (\ref{E-tauBnm}) reduces to the linear system
$\dt\tau_i +\frac{1}{2}\,\dx \tau_{i}=0$,
whose solution is a simple wave along $x$-axis: $\tau_i=\tau^0_i(t-2\,x)$.
Consider slightly more complicated cases.

1. For $f_j=\delta_{j2}$, (\ref{E-tauBnm}) is reduced to the system
\begin{equation}\label{E-tauk-1}
 \dt\tau_i = -\frac{i}{i+1}\,\dx\tau_{i+1},\quad i=1,2,\ldots
\end{equation}
The $n$-truncated (\ref{E-tauk-1}) reads as $\dt\vec\tau+B_{n,1}\,\dx\vec\tau=0$.
For $n=2$, using (\ref{E-symm-k1}), we have just two PDEs
\begin{equation*}
 \dt\tau_1 = -\frac{1}{2}\dx\tau_{2},\quad
 \dt\tau_2 = -\frac{2}{3}\dx\tau_{3}
 =(\tau_{1}^2-\tau_{2})\dx\tau_{1} - \tau_{1}\dx\tau_{2}.
\end{equation*}
The matrix
 $B_{2,1} =\left(\begin{array}{cc}
           0 & \frac12 \\
          -2\,\sigma_{2} & \sigma_{1} \\
           \end{array}\right)$
has the characteristic polynomial
$P_2=\det|\lambda \id_2-B_{2,1}|=\lambda^2-\sigma_1\lambda + \sigma_2$.
If $\lambda_1\ne\lambda_2$, the eigenvectors of $B_{2,1}$ are equal to $\vec v_j=(1,\, 2\lambda _j)\ (j=1,2)$.
If $\lambda_{1}=\lambda_{2}\ne0$, $B_{2,1}$ has one eigenvector only, hence the system (\ref{E-tauk-1})
is not hyperbolic in the $t$-direction.

For $n=3$, (\ref{E-tauk-1}) reduces to the quasilinear system of three PDEs with the matrix
$$B_{3,1}=\left(\begin{array}{ccc}
                     0 & \frac12 & 0 \\
                     0 & 0 & \frac23 \\
                     3\,\sigma_{3} & -\frac{3}{2}\,\sigma_{2} & \sigma_{1} \\
     \end{array}\right)
$$
 The characteristic polynomial of $B_{3,1}$ is
$ P_3=\lambda^3-\sigma_1\lambda^2+\sigma_2\lambda -\sigma_3,$
the eigenvalues are $\lambda_{j}$,
and the eigenvectors are $v_j=(1, 2\,\lambda_j, 3\,\lambda_j^2),\ j=1,2,3$.

\vskip1mm
2. For $f_j=\delta_{j3}$, (\ref{E-tauBnm}) reduces to the system
\begin{equation}\label{E-tauk-2}
 \dt\tau_i + \frac{3\,i}{2(i+2)}\,\dx\tau_{i+2}=0,\qquad i=1,2,\ldots
\end{equation}
The matrix of this $n$-truncated system is $B_{n,2}=\frac32(B_{n,1})^2$.

For $n=3$, (\ref{E-tauk-2}) reduces to the system of three quasilinear PDEs
\begin{equation*}
 \dt\tau_1 = -\frac{1}{2}\dx\tau_{3},\quad
 \dt\tau_2 = -\frac{3}{4}\dx\tau_{4}, \quad
 \dt\tau_3 =-\frac{9}{10}\dx\tau_{5},
 \end{equation*}
where $\dx\tau_{4}$  and $\dx\tau_{5}$ should be expressed using (\ref{E-symm-k1}) and (\ref{E-symm-k2}).
 The matrix of this system~is
\begin{eqnarray*}
 B_{3,2}\eq \frac32 (B_{3,1})^2 =\left(\begin{array}{ccc}
                    0 & 0 & \frac{1}{2} \\
                    3\,\sigma_3 &
                    -\frac{3}{2}\sigma_2 & \sigma_1 \\
                     \frac{9}{2}\sigma_1\sigma_3 &
                     \frac{9}{4}(\sigma_3-\sigma_1\sigma_2) &
                     \frac{3}{2}(\sigma_1^2-\sigma_2) \\
           \end{array}\right).
\end{eqnarray*}
It has eigenvalues $\tilde\lambda_j=\frac{3}2\lambda_j^{2}$, and the same eigenvectors as for $B_{3,1}$.

For $n=4$, (\ref{E-tauk-2}) reduces to the quasilinear system with the matrix
\begin{eqnarray*}
 B_{4,2}\eq\frac32 (B_{4,1})^2 =\left(\hskip-4pt\begin{array}{cccc}
           0 & 0 & \frac{1}{2} & 0\\
           0 & 0 & 0 & \frac{3}{4}\\
          -\frac92\,\sigma_4 & \frac94\,\sigma_3 & -\frac{3}{2}\sigma_2 & \frac98\sigma_1\\
           -6\,\sigma_1\sigma_4 &
           3\big(\sigma_1\sigma_3{-}\sigma_4\big) &
           2\big(\sigma_3{-}\sigma_1\sigma_2\big) &
           \frac32\big(\sigma_1^2{-}\sigma_2\big)\\
           \end{array}\hskip-4pt\right)
\end{eqnarray*}
with eigenvalues $\vec v_j=\frac{3}2\lambda_i^{2}$,
and eigenvectors $v_j=(1, 2\,\lambda_j, 3\,\lambda_j^2, 4\lambda_j^3),\ j=1,2,3,4$.


This series of examples can be continued as long as one desires.
\end{ex}

\subsection{Biregular foliated coordinates}

The~coordinate system described in the following lemma (for the proof see \cite{cc}, Section 5.1),
is called a \textit{biregular foliated chart} here.

\begin{lem}\label{L-adapted-coords}
Let $(M,\calf)$ be a differentiable manifold with a codimension-1 foliation
and a vector field $N$ transversal to $\calf$.
Then for any $p\in M$ there is a coordinate system $(x_0,x_1,\ldots x_n)$ on a neighborhood $U_p\subset M$ (centered at $p$) such that the leaves on $U_p$ are given by $\{x_0=c\}$ (hence the coordinate vector fields $\partial_i=\partial_{x_i},\ i\ge1$, are tangent to leaves), and $N$ is directed along $\partial_0=\partial_{x_0}$ (one may assume $N=\partial_0$ at $p$).
\end{lem}

If $(M, \calf , g)$ is a foliated Riemannian manifold and $N$ is a unit normal,
then $g$ has, in biregular foliated coordinates $(x_0, x_1, \ldots , x_n)$, the form
\begin{equation}\label{E-G0ij-g1}
 g=g_{00}\,dx_0^2+\sum\nolimits_{i,j>0}g_{ij}\,dx_i dx_j.
\end{equation}

\begin{lem}\label{lem-biregG}
For a metric (\ref{E-G0ij-g1}) in biregular foliated coordinates of a codimension one foliation $\calf$ on $(M,g)$, one has
 $N =\partial_0/\sqrt{g_{00}}$ (the unit normal)~and
\begin{eqnarray*}
 \Gamma ^{0}_{ij} \eq -\frac12\, g_{ij,0}/g_{00},\quad
 \Gamma ^j_{i0} = \frac12\sum\nolimits_{\/s} g_{is,0}g^{sj}
 \quad\mbox{(the Christoffel symbols)},\\
 b_{ij}\eq \Gamma^0_{ij}\sqrt{g_{00}}=-\frac12\,g_{ij,0}/\sqrt{g_{00}}
 \quad\mbox{(the second fundamental form)},\\
 A^j_{i}\eq -\Gamma^j_{i0}/\sqrt{g_{00}}=\frac{-1}{2\,\sqrt{g_{00}}}\sum\nolimits_{\/s} g_{is,0}\,g^{sj}
 \quad\mbox{(the Weingarten operator)},\\
 (b_m)_{ij}\eq
 A^{s_1}_{s_2}
 \dots A^{s_{m-1}}_{s_{m}} A^{s_{m}}_i g_{js_1}
 \quad\mbox{(the $m$-th ``power" of $\,b_{ij}$)},\\
 \tau_m\eq \Big(\frac{-1}{2\sqrt{g_{00}}}\Big)^m\sum\nolimits_{\,\{r_a\},\, \{s_b\}}
 g_{r_1 s_2,0}\ldots g_{r_{m-1} s_m,0}\,g^{s_1 r_1}\ldots g^{s_m r_m}.
\end{eqnarray*}
In particular,
\begin{equation*}
 \tau_1=\frac{-1}{2\sqrt{g_{00}}}\sum\limits_{\,r,\,s} g_{r s,0}\,g^{r s},\ \
 \tau_2=\frac{1}{4\,g_{00}}\sum\limits_{\,r_1,\, r_2,\, s_1,\, s_2}
 \hskip-3mm g_{r_1 s_2,0}\,g_{r_2 s_1,0}\,g^{s_1 r_1}\,g^{s_2 r_2},\ \ \mbox{\rm etc}.
\end{equation*}
\end{lem}

\noindent\textbf{Proof}.
All of that is standard and left to a reader. Just for his convenience, lest us observe that
the formula for $A$ follows from that for $b$ and $A^j_{i}{=}\sum\nolimits_{\/s}b_{is}g^{sj}$.
Notice that
 $(A^m)^j_i{=}\sum_{\{s_l\}}A^j_{s_2} A^{s_2}_{s_3}\dots
 A^{s_{m-1}}_{s_{m}} A^{s_{m}}_i$.
The formulae for $b_m$ follow from the above and
$(A^m)^s_i g_{sj}{=}g(A^m e_i,e_j){=}(b_m)_{ij}$.
The formulae for $\tau$'s follow directly from the above and the equality $\tau_m=\tr(A^m)$.\hfill$\square$

\section{A solution to general case}
\label{sec:gen-case}

\subsection{Searching for $\vec\tau$}
\label{sec:tau-general}

Let $g_t$ satisfy (\ref{eq1}) and $N_t$ be the $g_t$-unit normal
vector field to $\calf$ on $M$.

It is easy to see that $\dt N_t = 0$, therefore $N_t = N$ for all $t$, $N$ being the unit
normal of $\calf$ on $(M, g_0)$.
In fact, for any vector field $X$ tangent to $\calf$ one has
$0 = \dt g_t(X, N_t) = h(b_t)(X, N_t) + g_t(X, \dt N_t) = g_t(X, \dt N_t)$
and similarly
$0 = \dt g_t(N_t, N_t) = h(b_t)(N_t, N_t) + 2g_t(N_t, \dt N_t) = 2g_t(N_t, \dt N_t)$.

\vskip1mm
Let now $\nabla^t$ denote the Levi-Civita connection on $(M, g_t)$ and
$\Pi_t = \dt\nabla^t$. $\Pi_t$ is a (1,2)-tensor field on $M$.
Following, for example, \cite{to}, one can write for all $t$-independent $X,Y$ and $Z$
\begin{equation}\label{eq2}
 g_t(\Pi_t(X, Y), Z) =\frac{1}{2}\big[(\nabla^t_X h_t)(Y,Z)
 +(\nabla^t_Y h_t)(X,Z) -(\nabla^t_Z h_t)(X,Y)\big].
\end{equation}

\begin{lem}\label{P-Ak}
 For the EGF of type (b), on the tangent bundle of $\calf$ we have
\begin{eqnarray}\label{E-Ak}
\nonumber
 \dt b(X,Y) \eq h_t(A X, Y) \\
 &&\hskip-25mm-\frac12\sum\nolimits_{j=0}^{n-1}[
 \,N(f_j(\vec\tau,t))\,g_t(A^j X,Y)+f_j(\vec\tau,t)\,g_t((\nabla^t_N A^j)X,Y)],\\
 \label{E-Ak2}
 \dt A \eq-\frac12\sum\nolimits_{j=0}^{n-1}[N(f_j(\vec\tau,t))A^j +f_j(\vec\tau,t)\nabla^t_N A^j].
\end{eqnarray}
\end{lem}

\noindent\textbf{Proof}.
 By definition, $b(X,Y) = g_t(\nabla^t_X Y, N)$, and $h(\cdot, N)=0$.
 Using (\ref{eq2}) and the identity $h(A X,Y)=h(A Y,X)$, we obtain
\begin{eqnarray*}
 &&\dt b(X,Y) =
 \dt g_t(\nabla^t_X Y, N) = (\dt g_t)(\nabla^t_X Y, N) + g_t(\dt\nabla^t_X Y, N)\\
 &&=(1/2)\big[(\nabla^t_X h_t)(Y,N)+(\nabla^t_Y h_t)(X,N)-(\nabla^t_N h_t)(X,Y)\big] +h_t(\nabla^t_X Y, N)\\
 &&=(1/2)\big[h_t(A X,Y)+h_t(A Y,X)-N(h_t(X,Y))\big]
 =h_t(A X, Y) \\
 &&-(1/2)\sum\nolimits_{j=0}^{n-1}[f_j(\vec\tau,t)\,g_t((\nabla^t_N A^j) X,Y) +N(f_j(\vec\tau,t))\,g_t(A^j X,Y)].
\end{eqnarray*}
This proves (\ref{E-Ak}). Now, (\ref{E-Ak2}) follows from (\ref{E-Ak}) and
\begin{eqnarray*}
  g_t((\dt A)X, Y)\eq g_t(\dt(A X), Y)= \dt(g_t(A X, Y)) -(\dt g_t)(A X, Y)\\
 \eq \dt b(X, Y) -h_t(A X, Y).\hskip55mm \square
\end{eqnarray*}

The following lemma shows that the functions $\tau _1, \ldots , \tau _n$ satisfy the system of quasilinear PDEs,
whose matrix can be built using a generalized companion matrix of the characteristic polynomial of $A_N$.

\begin{lem}\label{P-tauAk}
The power sums $\{\tau_i\}_{i\in\NN}$ of the EGF (\ref{eq1}) of type (b) satisfy infinite system
of quasilinear PDEs (\ref{E-tauAk}).
 The $n$-truncated system (\ref{E-tauAk}) has the form
$\dt\vec\tau+(\tilde A+\tilde B)\,\dx\vec\tau=0$ with $\tilde{A}$ and $\tilde{B}$ given by (\ref{E-ABnm}).
\end{lem}

\noindent\textbf{Proof}. Applying $i A^{i-1}$ to both sides of (\ref{E-Ak2}), we obtain the PDE
\begin{equation*}
 i\,A^{i-1}\dt A = -\frac i2\sum\nolimits_{j=0}^{n-1}\big[N(f_j) A^{i+j-1}+f_j\,A^{i-1}\nabla^t_N A^{j}\big].
\end{equation*}
Taking the trace of both sides of the above equality, and using the identities
\begin{eqnarray*}
 i\tr(A^{i-1}\dt A) \eq \tr(\dt A^i)=\dt\tau_i,\\
 \frac{i+j-1}{j}\tr(A^{i-1}\nabla^t_N A^{j})\eq\tr(\nabla^t_N A^{i+j-1})=N(\tau_{i+j-1}),
\end{eqnarray*}
for $i,j>0$, we obtain (\ref{E-tauAk}).
The second statement of our lemma follows directly from Proposition \ref{L-BBnm}.
\hfill$\square$

\begin{rem}\rm
Lemma \ref{P-tauAk} implies the variational formula for $\sigma_m$.
Notice that for any differentiable matrix valued function $A(t)\in\mathrm{Mat}(n\times n)$, $t\in[0,\eps)$,
one has $\sigma_m'=\sum\nolimits_{i=0}^{m-1}(-1)^{i}\,\sigma_{m-i-1}\tr(A^{i} A')$,
$m=1,2,\ldots,n$, see \cite{rw2}.
In view of $(i+1)\tr(A^{i} A')=\tau_{i+1}'$, we obtain
\[
 \dt\sigma_m = \sum\nolimits_{i=0}^{m-1}(-1)^{i}\frac1{i+1}\,\sigma_{m-i-1}\dt\tau_{i+1}.
\]
 From above identity  with $\dt\tau_{i+j}$ replaced due to (\ref{E-tauAk}) follows the equality
\begin{eqnarray}\label{eq6}
\nonumber
\dt\sigma_m \eq\frac12\sum\nolimits_{i=0}^{m-1}{(-1)^{i+1}}\sigma_{m-i-1} \big\{ N(f_0(\vec{\tau},t))\,\tau_{i}\\
 \plus\sum\nolimits_{j=1}^{n-1}\big[ N(f_j(\vec{\tau},t))\,\tau_{i+j} +\frac{j}{i+j}f_j(\vec{\tau},t)\,N(\tau_{i+j})\big]\big\}.
\end{eqnarray}
\end{rem}

The next lemma deals with the evolution equation for EGF of type (a).

\begin{lem}\label{L-Tlocalsol2}
Let $g_t$ be the solution to the EGF (\ref{eq1}) of type (a).
Then the Weingarten operator $A$ of $\calf$ with respect to $g_t$ satisfies
\begin{equation}\label{E-Ak2-B}
 \dt A = -\frac12\big[N(\tilde f_0) \id +\sum\nolimits_{j=1}^{n-1}
 \big(N(\tilde f_j) A^j +\tilde f_j\cdot\nabla^t_N A^j\big)\big],
\end{equation}
and $\tau_i\ (i\ge1)$ (the power sums of the eigenvalues of $A$) satisfy PDEs
\begin{equation}\label{E-tauAk-tilde}
 \dt\tau_i=-\frac i2\big\{\tau_{i-1}N(\tilde f_0) +\sum\nolimits_{j=1}^{n-1}
 \big[\,\frac{j}{i+j-1}\tilde f_j\,N(\tau_{i+j-1}) +\tau_{i+j-1}N(\tilde f_j)\big]\big\} .
\end{equation}
The $n$-truncated system (\ref{E-tauAk-tilde}) has the form
\begin{equation*}
 \dt\vec\tau+\big[\sum\nolimits_{j=1}^{n-1} j\,\tilde f_j\,(B_{n,1})^{j-1}\big] N(\vec\tau)=\vec{a},
\end{equation*}
where $\vec{a} = (\hat a_1, \ldots , \hat a_n)$,
$\hat a_i=-\frac{i}2\sum\nolimits_{j=0}^{n-1} N(\tilde f_j)\,\tau_{i+j-1}\ (1\le i\le n)$,
and $B_{n,1}$ is the generalized companion matrix (\ref{E-Bn1}).
\end{lem}

\noindent\textbf{Proof}. The proof of (\ref{E-Ak2-B}) is similar to that of (\ref{E-Ak2}).
On $T\calf$ we obtain
\begin{eqnarray*}
 g_t((\dt A)X, Y)\eq-\frac12 N(h_t(X,Y)) =-\frac12\big[N(\tilde f_0)g_t(X,Y)\\
 &&+\sum\nolimits_{j=1}^{n-1}
 \big(N(\tilde f_j)g_t(A^j X,Y) +\tilde f_j\,g_t(\nabla^t_N A^jX,Y)\big)\big].
\end{eqnarray*}
Following the lines of the proof of (\ref{E-tauAk}) in Lemma~\ref{P-tauAk},
we deduce from above the system~(\ref{E-tauAk-tilde}).
By Proposition~\ref{L-BBnm}, the $n$-truncated (\ref{E-tauAk-tilde}) has the required form.$\hfill\square$

\vskip1.5mm
Lemmas \ref{P-tauAk} and \ref{L-Tlocalsol2} together with Theorem A provide existence and uniqueness results for the symmetric functions $\vec{\tau}^{\,t}$ satisfying conditions following from (\ref{eq1}). In particular, this allows to reduce the problems of existence and uniqueness for EGF of type (b) to those for type (a) as we do in the proof of Theorem \ref{T-localsol3} in the next section.

\subsection{Local existence of $g_t$ (Proofs of Theorems~\ref{L-main1} and \ref{T-localsol3})}
\label{sec:Texist}

Given a Riemannian metric $g$ on a foliated manifold $(M,\calf)$, the symmetric tensor $h(b)$ defined by (\ref{eq1}) can be expressed locally in terms of the first partial derivatives of~$g$.
Therefore, $g\mapsto h(b)$ is a first order partial differential operator.
For a $\calf$-truncated symmetric $(0,2)$-tensor $\hat S$ on $T\calf$,
\begin{equation}\label{E-Acoord-0}
 h(b)=\hat S
\end{equation}
is a non-linear system of first order PDEs.
 The particular case of (\ref{E-Acoord-0}) is the Einstein type relation
 $h(b)=\lambda\,\hat g$ for some function (or constant) $\lambda$ on $M$.

\vskip1.5mm
There are several obvious obstructions to the existence of solutions to (\ref{E-Acoord-0}) even at a single point:
for example, if $h(b)=\hat b_{2j}$
and $\hat S$ takes both, positive and negative, values at $p\in M$, then (\ref{E-Acoord-0}) has no solutions at $p$.

\vskip1.5mm
\noindent\textbf{Proof of Theorem~\ref{L-main1}}.
Take biregular foliated coordinates\index{biregular foliated coordinates} $(x_0,x_1,\ldots x_n)$ on $U_q\subset M$ (with center  at $q$); see Lemma~\ref{L-adapted-coords}, and the metric (\ref{E-G0ij-g1}). Then, $N=\partial_0/\sqrt{g_{00}}$ is the unit normal to $\calf$.
Set $\psi_{ab}=g_{ab,0}$.
The system (\ref{eq1}) (for $f_j$ of type (a)) along a trajectory
$\gamma: x\mapsto \gamma(x)$ of $\partial_0$ has the form
\begin{equation}\label{E-local2-1}
 \dt g_{ij}=F_{ij}(g_{ab}, \psi_{ab}, t, x),
\end{equation}
where $F_{ij}:= h(b)_{ij}$.
In view of symmetry, $\psi_{ab}=\psi_{ba}$ and $F_{ij}=F_{j\/i}$, we shall assume $1\le i\le j\le n$ and $1\le a\le b\le n$.
For example, if $f_m=0\ (m>1)$ then (\ref{E-local2-1}) is the hyperbolic (diagonal) system
\begin{equation*}
 \dt g_{ij} = f_0(q,t) g_{ij}-\frac12 g_{00}^{-1/2}f_1(q,t) \psi_{ij},
\end{equation*}
that completes the proof in this case.

Now let $f_m\ne0$ for some $m>1$ (e.g., general $f_m$).
We may assume that
$A\partial_j=k_j\partial_j$, $g(\partial_i,\partial_j)=\delta_{ij}$ ($i,j>0$) and $g_{00} =1$
at the point $q$ for $t=0$.
(By Lemma \ref{lem-biregG}, we have
 $(b_m)_{ij}=
 (-1/2
 )^m
 \psi^m_{ij}\,\delta_{ij}$
at $q$ for $t=0$).

Differentiating (\ref{E-local2-1}) with respect to $x$ and $t$, we obtain
\begin{eqnarray}\label{E-pF1}
\nonumber
 \partial_0 p_{ij} \eq\partial_0 F_{ij}
 +\sum\nolimits_{a,b}\Big[\frac{\partial F_{ij}}{\partial g_{ab}}\,\partial_0 g_{ab}
 +\frac{\partial F_{ij}}{\partial \psi_{ab}}\,\partial_0\psi_{ab}\Big],\\
 \dt p_{ij} \eq\dt F_{ij}
 +\sum\nolimits_{a,b}\Big[\frac{\partial F_{ij}}{\partial g_{ab}}\,\dt g_{ab}
 +\frac{\partial F_{ij}}{\partial \psi_{ab}}\,\dt\psi_{ab}\Big],
\end{eqnarray}
where
$ p_{ij}:=\dt g_{ij},\
 \psi_{ij}:=g_{ij,0},\
 F_{ij}:= h(b)_{ij}$.

As  $g$ is of  class $C^2$, we conclude that
 $\dt\psi_{ab}=\frac{\partial^2 g_{ab}}{\partial t\partial x_0}=\partial_0 p_{ab}$.
Hence, (\ref{E-pF1}) together with (\ref{E-local2-1}) may be written in the~form
\begin{eqnarray}\label{E-pF2}
\nonumber
 \dt g_{ij}\eq F_{ij}(\{g_{ab}\}, \{\psi_{ab}\}, t, x),\\
 \dt\psi_{ij} -\sum\nolimits_{a,b}\frac{\partial F_{ij}}{\partial \psi_{ab}}\,\partial_0\psi_{ab}
 \eq\partial_0 F_{ij}+\sum\nolimits_{a,b}\frac{\partial F_{ij}}{\partial g_{ab}}\,\psi_{ab},\\
\nonumber
 \dt p_{ij} -\sum\nolimits_{a,b}\frac{\partial F_{ij}}{\partial\psi_{ab}}\,\partial_0 p_{ab}
 \eq\dt F_{ij}+\sum\nolimits_{a,b}\frac{\partial F_{ij}}{\partial g_{ab}}\,p_{ab}.
\end{eqnarray}
The above quasilinear system consists of parts:

(i) our original equation (\ref{E-pF2})$_1$,

(ii) the corresponding equation (\ref{E-pF2})$_2$ for $\dt A$,  and

(iii) the equation (\ref{E-pF2})$_3$ for $\dt^2 g$ following from the previous ones.

\noindent
In general, the following $\frac12n(n+1)\times\frac12n(n+1)$ matrix is not symmetric:
\[
 d_{\psi}F=\Big\{\frac{\partial F_{ij}}{\partial\psi_{ab}}\Big\},\quad
 i\le j,\ \ a\le b.
\]
We claim that it is hyperbolic.
If we change the local coordinate system on $M$, the components $F_{ij}\ (i\le j)$ and $\psi_{ab}\ (a\le b)$ at $q$ will be transformed by the same tensor low.
Notice that the above $d_{\psi}F$ is a $(1,1)$-tensor on the vector bundle of symmetric $(0,2)$-tensors on $T\calf$.
Hence,  $d_{\psi}F(q)$
can be seen as the linear endomorphism of the space of symmetric $(0,2)$-tensors on $T_q\calf$.

The hyperbolicity is a pointwise property, so can be considered at any point $q\in
M$ in a special bifoliated chart around $q$,
for example, such that $g_{ij} = \delta_{ij}$ and $A_{ij}(q) = k_i\delta_{ij}$
at $q$.
(Indeed, $(k_i)$ are the principal curvatures of $\calf$ at $q$ for $t=0$).
In this chart, our calculations show that the matrix $d_{\psi}F$ is diagonal,
so has real eigenvalues (vectors) at $q$.
Indeed, for $t=0$ one may find at the point $q$:
$$
 \frac{\partial F_{ij}}{\partial\psi_{ab}}=
 \sum\limits_{m\ge1} f_m(q)\,\mu(m)_{ij}\delta^{\{i,j\}}_{\{a,b\}}, \qquad
 \mu(m)_{ij}=\sum\limits_{\alpha+\beta=m-1} k_i^\alpha k_j^\beta.
$$
The order of indices of $d_{\psi}F$ is $[1,1],[1,2],\dots,[1,n],[2,2],
\dots,[2,n],\dots,[n,n]$.
For example, for $F=b_2$ (i.e., $f_j=\delta_{j2}$) the above matrix in an orthonormal frame
at any point is
\[
\frac{\partial (b_2)_{ij}}{\partial\psi_{ab}}=
\frac14
\left[\begin{smallmatrix} 2\,\psi_{11}&2\,\psi_{12}&2\,\psi_{13}&0&0&0\\ \noalign{\medskip}\psi_{12}&\psi_{11}+\,\psi_{22}&\psi_{23}&\psi_{{12}}&\psi_{13}&0\\ \noalign{\medskip}\psi_{13}&\psi_{23}&\psi_{11}+\,\psi_{33}&0&\psi_{12}&\psi_{13}\\ \noalign{\medskip}0&2\,\psi_{12}&0&2\,\psi_{22}&2\,\psi_{23}&0\\ \noalign{\medskip}0&\psi_{13}&\psi_{{12}}&\psi_{23}&\psi_{22}+\,\psi_{33}&\psi_{23}\\
\noalign{\medskip}0&0&2\,\psi_{13}&0&2\,\psi_{23}&2\,\psi_{33}
\end{smallmatrix}\right]
\]
with the order of indices $[1,1],[1,2],[1,3],[2,2],[2,3],[3,3]$.
At $q$ for $t=0$ (i.e., $\psi_{ab}=0,\ a\ne b$ and $\psi_{aa}=k_a$) it is diagonal with the elements $\mu(2)_{ab}=k_a+k_b$.
 Let $A_0=[\frac12, 1,1, \frac12, 1,\frac12]$ be the diagonal matrix.
 Then the matrix $A_1=A_0\,d_\psi (b_2)$ is symmetric, and our system (\ref{E-pF2})$_2$ for $h(b)=b_2$ is ``symmetrizable":
 $A_0\dt\psi - A_1\partial_0\psi =\mbox{\{free terms\}}$.

Therefore, (\ref{E-pF2}) for the functions $g_{ij}(t,x),\,p_{ij}(t,x)$ and $\psi_{ij}(t,x)$ with $i\le j$ is the hyperbolic system,
which is ``symmetrizable".
Indeed, multiplying $n$ columns  (corresponding to $i=j$)
of the matrix $d_\psi F$ in an orthonormal frame by the factor $\frac12$,  we obtain the symmetric matrix.

 By Theorem~A (in Section~\ref{sec:quasi}),  given $q\in M$ there exists a unique solution to (\ref{E-pF2}) which is defined in $U_q$ along the $N$-curve through $q$ for some time interval $[0,\eps_q)$ and satisfies the initial conditions
 \[g_{ij}(0,x)=(g_{0})_{ij}(x),\quad p_{ij}(0,x)=h(b_0)_{ij}(x),\quad
 \psi_{ij}(0,x)=(\partial_0 g)_{ij}(x).
\]
By Theorem~A, see definitions of $K$ and $\bar\Omega$, the value $\eps_q$ continuously depends on $q\in M$.
 The claim follows from the above and compactness of $M$.
 $\hfill\square$

\vskip1.5mm
\noindent\textbf{Proof of Theorem~\ref{T-localsol3}}.
Let $A_0$ and $\vec\tau^{\,0}$ be the values
 of extended Weingarten operator and power sums of the principal curvatures $k_i$ of (the leaves of) $\calf$
 determined on $(M,\calf)$ by a given metric $g_0$.

 a) {\it Uniqueness}. Assume that $g^{(1)}_t,g^{(2)}_t$ are two solutions to (\ref{eq1}) with the same initial metric $g_0$.
 The functions $\vec\tau^{\,t,1},\vec\tau^{\,t,2}$, corresponding to $g^{(1)}_t,g^{(2)}_t$, satisfy (\ref{E-tauAk}) and have the same initial value $\vec\tau^{\,0}$.
 By~Lemma~\ref{P-tauAk} and Theorem~A, $\vec\tau^{\,t,1}=\vec\tau^{\,t,2}_t=\vec\tau^{\,t}$ on some positive time interval $[0,\eps_1)$. Hence $g^{(1)}_t,g^{(2)}_t$ satisfy (\ref{eq1}) of type (a) with known $\tilde f_j(p,t):=f_j(\vec\tau^{\,t}(p),t)$. By Theorem~\ref{L-main1}, $g^{(1)}_t=g^{(2)}_t$ on some positive time interval $[0,\eps_2)$.

 b) {\it Existence}.
  By Lemma~\ref{P-tauAk} and Theorem A, there is a unique solution $\vec\tau^{\,t}$ to (\ref{E-tauAk}) on some positive time interval $[0,\eps^*)$.
 By Theorem~\ref{L-main1}, the EGF (\ref{eq1})  of type (a) with known functions $\tilde f_j(\cdot,t):=f_j(\vec\tau^{\,t},t)$ has a unique solution $g^*_t$ ($g^*_0=g_0$) for $0\le t<\eps^*$.
 The Weingarten operator $A^*_t\ (A^*_0=A_0)$ of $(M,\calf,g^*_t)$ satisfies (\ref{E-Ak2-B}),
 hence the power sums of its eigenvalues, $\vec\tau^{\,t,*}\ (\vec\tau^{\,0,*}=\vec\tau^{\,0})$, satisfy (\ref{E-tauAk-tilde}) with the same coefficient functions $\tilde{f}_j$.
 By Lemma~\ref{L-Tlocalsol2} and Theorem A, the solution of this problem is unique,
 hence $\vec\tau^{\,t}=\vec\tau^{\,t,*}$,
 i.e., $\vec\tau^{\,t}$ are power sums of eigenvalues of $A^*_t$.
 Finally, $g^*_t$ is a solution to (\ref{eq1}) such that
 $\vec\tau^{\,t}$ are power sums of the principal curvatures of the leaves
 in this metric.$\hfill\square$

\subsection{Proofs of Corollaries~\ref{C-A1}, \ref{C-A0} and Propositions~\ref{T-f1-B}\,--\,\ref{C-tu1}}

In all the proofs below, $p$ is an arbitrary point of $M$
and $\gamma : x\mapsto \gamma(x)$ $(\gamma(0)=p,\ x\in\RR)$ is the trajectory of $N$, $N$ -- the unit normal of $\calf$.

\vskip1.5mm
\noindent\textbf{Proof of Corollary~\ref{C-A1}}.
By Lemma~\ref{P-tauAk}, we have (\ref{E-tauAk}), which in our case  reduces to
the system~(\ref{eq322}). The last one means the initial value problem in the $(x,t)$-plane
for the vector function $\vec\tau(x,t)=\vec\tau(\gamma(x),t)$
\begin{equation}\label{E-PDE-1b}
 \dt\vec\tau + \big(\frac12 f(\vec\tau,t)\id_n +\tilde A(\vec\tau,t)\big)\dx\vec\tau=0,\quad
 \vec\tau(x,0)=\vec\tau\,^0(\gamma(x)).
\end{equation}
The matrix $\tilde A$ is equal to
$\{\frac i2\,\tau_{i} f_{,\tau_j}(\vec\tau,t)\}_{1\le i,j\le n}$. Note that ${\rm rank}\,\tilde A \le 1$.
By condition ($H_2$), either the function
$\tilde\lambda=\tr\tilde A=\sum\nolimits_{1\le i\le n}\frac i2\,\tau_i f_{,\tau_i}(\vec\tau,0)$
(the eigenvalue of~$\tilde A$) is non-zero for all $x\in\RR$,
or $\tilde A(x)\equiv0$. Hence (\ref{E-PDE-1b}) is hyperbolic for small enough $t$.
In first case, the eigenvector of $\tilde A_{|\,t=0}$ for $\tilde\lambda(x)$
is $\bv_1=(f_{,\tau_1}, f_{,\tau_2},\ldots, f_{,\tau_n})$, and the kernel of $\tilde A_{|\,t=0}$
is spanned by $n-1$ vectors
$$
 \bv_2{=}(-2\,\tau_{2}, \tau_{1}, 0,\ldots 0),\
 \bv_3{=}(-3\,\tau_{3}, 0, \tau_{1}, 0,\ldots 0),\ldots
 \bv_n{=}(-n\,\tau_{n}, 0,\ldots 0, \tau_{1}).
$$
(If $\tilde\lambda(x)=0$ for some $x\in\RR$, and $f_{,\tau_j}(\vec\tau,0)\ne0$ for some $j$,
then $\tilde A$ is nilpotent and hence (\ref{E-PDE-1b}) is not hyperbolic).
 By Theorem~A, the initial value problem (\ref{E-PDE-1b}) has a unique solution on a domain
$[-\delta,\delta]\times[0,\eps')$ of the $(x,t)$-plane. Hence there exists $t_p > 0$ such that
the solution $\vec\tau(\cdot,t)$ to (\ref{eq322}) exists and is unique for $t\in[0,t_p)$ on a neighborhood  $U_p\subset M$ centered at $p$. By compactness of $M$, we conclude that there is $\eps>0$ such that
(\ref{eq-g-f1}) admits a unique solution $\vec\tau(p,t)$ on $M$ for $t\in[0,\eps)$.

\vskip1mm
One may also apply the \textit{method of characteristics} to solve (\ref{E-PDE-1b}) explicitly when $f=f(\vec\tau)$.
For $\tilde\lambda\ne0$, the system has two characteristics:
$\frac{d\tilde x}{dt}=\tilde\lambda+\frac f2$ and $\frac{d x}{dt}=\frac f2$
corresponding to two eigenvalues $\tilde\lambda+\frac f2$ (of multiplicity 1) and $\frac f2$ (of multiplicity $n-1$).
 First, let us observe that the function $u:=\bv_1^{\,T}\cdot\vec\tau=\sum_{j\le n}f_{,\tau_j}\tau_j$
is constant along the first family of characteristics:
\begin{equation}\label{E-first}
 \frac{d}{dt} u= \dt u+(\tilde\lambda+\frac f2)N(u)=0 \Leftrightarrow
 u=const \ {\rm along} \
 \frac{d}{dt}\tilde x=\tilde\lambda+\frac f2.
\end{equation}

Next, let us find a function that is constant along the second family of characteristics.
For each $m>1$, due to the form of $\bv_m$, calculate the sum of the first equation in (\ref{E-PDE-1b}) multiplied by $-m\,\tau_m$ with the $m$-th equation multiplied by $\tau_1$
(along the trajectories of $\frac{dx}{dt}=\frac f2$)
\begin{equation*}
 \tau_1\dt\tau_{m} - m\tau_m\dt\tau_{1}
 +\frac{f}2\big(\tau_1N(\tau_{m}) - m\tau_m N(\tau_{1})\big)=0.
\end{equation*}
(The terms with $\sum_j f_{\tau_j} N(\tau_j)$ cancel.)
Using
$\frac{d}{dt}\tau_m=\dt\tau_{m}+N(\tau_{m})\frac{dx}{dt}=\dt\tau_{m}+\frac12 f N(\tau_{m})$,
we get (again, along the second family of characteristics)
\begin{equation}\label{E-first3}
 \tau_1\frac{d}{dt}\tau_m -m\,\tau_m \frac{d}{dt}\tau_1=0 \quad(m\ge2).
\end{equation}
The complete integral of (\ref{E-first3}) is
\begin{equation}\label{E-first4}
 \tau_m = C_m(x)\,\tau_1^m \quad(m\ge2),
\end{equation}
where $C_m(x)=\tau_m(x,0)/\tau_1^m(x,0)$ are known functions.
Since the smooth functions $\vec\tau(x,t)$ and $f(\vec\tau,t)$ exist for $t\in[0,t_p)$,
the EGF under consideration exists and is unique on $M\times[0, t_p)$.

\vskip1mm
Using the basic properties of the Lie derivative of $g_t$ along $N$,
one may show that (\ref{eq-g-f1}) is equivalent to
$\mathcal{L}_{Z_t}\, g_t =0$ with $Z_t=\dt+\frac12f(\vec\tau\,^t,t)\,N$.
Denote by $\Phi_t$ the flow of $\frac12 f(\vec\tau,t)N$. The solution $g_t$
can be defined by $g_t(\Phi_t X, \Phi_t Y) = g_0(X, Y)$, $g_t(X, N)=0$ and $g_t(N,N)=1$
for all $X$ and $Y$ tangent to~$\calf$. The solution
exists and is unique as long as the solution to (\ref{eq322}) does.
If~$f(\cdot,t)=$ const, the solution exists (and is unique) for all $t\ge 0$.$\hfill\square$
\vskip1.5mm

\noindent\textbf{Proof of Corollary~\ref{C-A0}}.
Denote by $\vec\tau(x,t)=(\tau^t_1,\ldots,\tau^t_n)$ the power sums of the principal curvatures of $\calf$ at the point $\gamma(x)$ in time $t$.
By Lemma~\ref{P-tauAk}, we have (\ref{E-tauAk}), which in our case reduces to
(\ref{E-PDE-1b-0a}). The last one provides the initial value problem in the $(x,t)$-plane
\begin{equation}\label{E-PDE-1b-0}
 \dt\vec\tau + \tilde A(\vec\tau,t)\,\dx\vec\tau=0, \quad
 \vec\tau(x,0)=\vec\tau\,^0(\gamma(x)),
\end{equation}
where $\tilde A$ is equal to
$\{\frac i2\,\tau_{i-1} f_{,\tau_j}(\vec\tau,t)\}_{1\le i,j\le n}$. As before, ${\rm rank}\,\tilde A \le 1$.
 Consider the function
 $\tilde\lambda=\tr\tilde A=\sum\nolimits_{1\le i\le n}\frac i2\,\tau_{i-1}f_{,\tau_i}(\vec\tau,0)$
 (the eigenvalue of~$\tilde A$).
By condition $(H_1)$, either $\tilde\lambda(x)\ne0$ for all $x\in\RR$,
 or $\tilde A(x)\equiv0$. Hence the system (\ref{E-PDE-1b-0}) is hyperbolic at $(x,0)$.
 As in the proof of Corollary~\ref{C-A1} we conclude that there is $\eps>0$ such that $\vec\tau(p,t)$ on $M$ exists and is unique for $t\in[0,\eps)$.
Certainly, a unique smooth solution to (\ref{E-gC-A0}) has the required form.$\hfill\square$

\vskip1.5mm
\noindent\textbf{Proof of Proposition~\ref{T-f1-B}}.
Notice that
$$
 \sum\nolimits_{j}(\vec\tau^{\,\bm\alpha})_{,\tau_j}\tau_j=
 m\,\vec\tau^{\,\bm\alpha},\quad
 \sum\nolimits_{j}j(\vec\tau^{\,\bm\alpha})_{,\tau_j}\tau_j=
 l\,\vec\tau^{\,\bm\alpha}.
$$
If~$f=\sum\limits_{{\bm\alpha}}c_{\bm\alpha}\,\vec\tau^{\,\bm\alpha}$
then $N(f)=\sum\limits_{j}\sum\limits_{{\bm\alpha}}c_{\bm\alpha}\vec\tau^{\,\bm\alpha}_{,\tau_j}N(\tau_j)$ (the derivative of $f$ along $N$).
\newline
 Define $\tilde f(x,t)=f_t(\vec\tau(\gamma(x)))$ and $\tilde f_0=\tilde f(\cdot,0)$.
 One has PDEs (\ref{E-PDE-1b}) in the $(x,t)$-plane.
The characteristics of the first family, see (\ref{E-first}), are lines,
and $\tilde f=const$ along them.
 To~show this, observe that (by definition of $f_t(\vec\tau)$)
$$u:=\sum\nolimits_{j} f_{t,\tau_j}(\vec\tau)\,\tau_j=m\tilde f(x,t),\quad
 \tilde\lambda:=\sum\nolimits_{j}\frac j2\,f_{t,\tau_j}(\vec\tau)\tau_j
 =\frac{l}2\,\tilde f(x,t).
$$
Since $\tilde f=u/m$ is constant along the first family of characteristics in the $(x,t)$-plane,
these characteristics are lines given by the equation
\begin{equation*}
 \frac{d}{dt}x=\frac12(l+1)\tilde f\quad\Leftrightarrow\quad x = \xi + \frac12(l+1)\tilde f_0(\xi)\,t.
\end{equation*}
Notice that $\tilde f\,'= \sum\nolimits_{\,{\bm\alpha}\in J_{m,l}}c_{\bm\alpha}\sum_{j} (\vec\tau^{\,0})^{\,\bm\alpha}_{,\tau_j}N(\tau^0_j)$.
If~$\tilde f_0\/'>0$ on $\gamma$, the solution $\tilde f$
exists for all $t\ge0$ (see Example~\ref{Ex-burgers}) and we set $t_p=\infty$.
If~$\tilde f_0\/'$ is negative somewhere on $\gamma$, then $\tilde f$ exists (and is continuous) for
$t\in[0,t_p)$ where $t_p=-2/[(l+1)\min_{x}\tilde f\,'_0(x)]$.

The second family of characteristics, $\frac{d}{dt}x=\frac f2$, also exists for $t\in[0,t_p)$.
To~show this, assume the opposite, i.e., there are $t_0\in (0,t_p)$ and a
trajectory $\gamma_1(t)$ of the second family of characteristics that cannot be continued for
values $t\ge t_0$. Therefore, the inclination $\tilde f(\gamma_1(t),t)/2$ of
$\gamma_1$ approaches to infinity when $t\to t_0$, a contradiction to
continuity of $\tilde f$ on the strip $t\in[0,t_0]$ in the $(x,t)$-plane.
$\hfill\square$

\vskip1.5mm
\noindent\textbf{Proof of Proposition~\ref{P-08}}.
For the EGF (\ref{eq1}) define the function $\psi$~by
\begin{equation}\label{E-psi-f}
 \psi(\lambda,t):=\sum\nolimits_{j=0}^{n-1} f_j(n\lambda,n\lambda^2,\ldots,n\lambda^n;t)\,\lambda^j.
\end{equation}
Since $\calf$ is totally umbilical for $t=0$,
$A_0=\lambda_0\,\widehat{\id}$.
Assuming $A=\lambda_t\,\widehat{\id}$ for small enough $t$,
we see that (\ref{E-Ak2}), or equivalently, (\ref{E-tauAk}) for $i=1$, yields the PDE
$\dt\lambda_t+\frac12\,N(\psi(\lambda_t,t))=0$, whose unique solution exists for small enough~$t\ge0$.
From the uniqueness of solution to (\ref{E-Ak2}) it follows that
$A=\lambda_t\,\widehat{\id}$.$\hfill\,\square$

\vskip1.5mm
\noindent\textbf{Proof of Proposition~\ref{C-tu1}}.
Recall that $A X=\lambda_t X$. Theorem~A (with $n=1$) provides the short-time existence and uniqueness of the solution $\lambda_t$ to (\ref{E-tauAk-L}) for $0\le t<T$.
Furthermore, the EGF is expressed as $\dt g_t=\psi(\lambda_t)\hat g_t$,
and the solution $\hat g_t$ has the required form.

Consider the function $\tilde\lambda(x,t)=\lambda(\gamma(x),t)$ in the $(x,t)$-plane
along the trajectory $\gamma(x),\ \gamma(0)=p$, of $N$,
and set $\tilde\lambda_0(x)=\lambda(\gamma(x),0)$.
The equation (\ref{E-tauAk-L}) in this case has the form of a \textit{conservation law}
\begin{equation}\label{E-tauAk-K}
 \dt\tilde\lambda +\frac12\,\dx(\psi(\tilde\lambda)) = 0.
\end{equation}
One may show the following;
\textit{If $\psi'(\tilde\lambda_0(x)),\,\tilde\lambda_0(x)\in C^1(\RR)$ and if $\tilde\lambda_0(x)$ and $\psi'(\tilde\lambda_0(x))$ are either non-decreasing or  non-increasing, the problem
\begin{equation}\label{E-tauAk-Kxt}
 \dt\tilde\lambda +\frac12\,\dx(\psi(\tilde\lambda)) = 0,\quad
 \tilde\lambda(x,0)=\tilde\lambda_0(x),\quad t\ge0
\end{equation}
has a unique smooth solution defined implicitly by the parametric equations},
\begin{equation}\label{E-cu}
 \tilde\lambda(x,t)=\tilde\lambda_0(x),\quad
 x=\xi+\frac12\,\psi'(\tilde\lambda_0(\xi))\,t.
\end{equation}
\textit{If~$\frac{d}{dx}\,\psi'(\tilde\lambda_0(x))$ is negative elsewhere along $\gamma$, then $\tilde\lambda(x,t)$ exists for}
\begin{equation*}
 t<t_p=-2/\inf\limits_{x\in\RR}\frac{d}{dx}\,\psi'(\tilde\lambda_0(x)).
\end{equation*}
For $\psi(\lambda)=\lambda^2$, (\ref{E-tauAk-K}) is reduced to the Burgers' equation,
see Section~\ref{sec:quasi}.$\hfill\square$

\section{Applications and Examples}
\label{sec:appl}

\subsection{Extrinsic Ricci  and Newton transformation flows}
\label{sec:extRic}

\vskip1.5mm
The \textit{extrinsic Riemannian curvature tensor} $\mbox{Rm}^{\rm\,ex}$ of $\calf$
is, roughly speaking, the difference of the curvature tensors of $M$ and of the leaves.
More precisely, by the Gauss formula, we have
$$\mbox{Rm}^{\rm\,ex}(Z, X)Y = g(A X,Y)A Z - g(A Z, Y)A X.$$
In Section \ref{sec:extRic}, we study the \textit{extrinsic Ricci flow}
for small dimensions $n>1$,
\begin{equation}\label{E-extRic}
\dt g_t=-2\Ric^{\rm\,ex}(b_t).
\end{equation}
The \textit{extrinsic Ricci} tensor is given by
$\Ric^{\rm\,ex}(X,Y)=\tr\mbox{Rm}^{\rm\,ex}(\cdot,X)Y$, where $X,Y\in T\calf$.
Hence,
\begin{equation}\label{E-ricA2}
  \Ric^{\rm\,ex}(b) = \tau_1\hat b_1 - \hat b_2,\quad
  (\Ric^{\rm\,ex}(b)=\sigma_2\,\hat g \ \ \mbox{ when } \ n=2).
\end{equation}
Therefore, $-2\Ric^{\rm\,ex}(b)$ relates to $h(b)$ of (\ref{eq1}) with $f_1=-2\tau_1,\,f_2=2$ (others $f_j=0$).

\begin{lem}
The principal curvatures of extrinsic Ricci flow satisfy PDEs
\begin{equation}\label{E-ExtRick}
  \dt k_i = N(k_i(\tau_1 - k_i)),\quad i=1,\ldots, n
\end{equation}
which for a totally umbilical $\calf$, i.e., $k_i=\lambda$, are reduced to the PDE
\begin{equation}\label{E-ExtRickumb}
  \dt\lambda = 2\,(n-1)\,\lambda\,N(\lambda).
\end{equation}
The extrinsic scalar curvature of the flow (\ref{E-extRic}) is $-2\,\sigma_2$, where
\begin{equation}\label{E-Ak2-Rex2}
 \dt\sigma_2=(\tau_{1}^2+\tau_{2})\,N(\tau_1)-\tau_1\,N(\tau_2)-\frac{2}{3}\,N(\tau_{3}).
\end{equation}
For $n=2$, (\ref{E-Ak2-Rex2}) reads as
$\dt\sigma_2=\tau_1\,\dx\sigma_2$.
\end{lem}

\noindent\textbf{Proof}.
By Lemma~\ref{P-tauAk}, for the extrinsic Ricci flow (\ref{E-extRic}) we have
\begin{equation}\label{E-Rictau}
 \dt\tau_i=i\,\tau_i N(\tau_1)+\tau_1N(\tau_i)-\frac{2\,i}{i+1}N(\tau_{i+1}),\quad
 i>0.
\end{equation}
Replacing $\tau_i=\sum_{j=1}^{\,n} (k_j)^i$ in (\ref{E-Rictau}), we obtain (\ref{E-ExtRick}).
Differentiating identity $2\sigma_2=\tau_1^2-\tau_2$, we obtain
 $2\,\dt \sigma_2=2\,\tau_1\,\dt\tau_1-\dt\tau_2$.
From this and (\ref{E-Rictau}) with $i=1,2$ it follows (\ref{E-Ak2-Rex2}).
\hfill$\square$

\begin{cor}\label{C-ExtRic1}
Let $(M, g_0)$ be a compact Riemannian manifold with a codi\-men\-sion-1 foliation
$\calf$ and a unit normal $N$.
Then there exists a unique solution $g_t,\ t\in[0,\eps)$
(for some $\eps>0$), to the extrinsic Ricci flow (\ref{E-extRic}) in any of the following cases (i) and (ii):

\vskip1mm
 (i) $n=2$, and \ $\tau_{1}\ne0$;\quad
 (in this case, $\hat g_t= \hat g_0\exp\big(\int_0^t \sigma_2\,dt\big)$);

\vskip1mm
 (ii) $n=3$, and \ $|\sigma_1|^3>27|\sigma_{3}|>0$.
\end{cor}

Notice that any, either positive or negative definite, operator $A_N$
satisfies inequalities Corollary~\ref{C-ExtRic1} (ii) if only $A_N$ is not proportional to the identity.

\vskip1.5mm
\noindent\textbf{Proof}.
It is sufficient to show that in conditions of theorem, (\ref{E-Rictau}) is hyperbolic in the $t$-direction.

\underline{Let $n=2$}. By equality $\tau_3=\frac32\tau_1\tau_2{-}\frac12\tau_1^3$,
the matrix of 2-truncated system (\ref{E-Rictau}) is
$
 C_2 =\Big(\begin{matrix}
                    -2\tau_{1} & 1        \\
                    -2\tau_{1}^2 & \tau_{1} \\
           \end{matrix}\Big)
$.
One can see directly, or applying condition $(H_1)$ to $\Ric^{\rm\,ex}(b)=\sigma_2\,\hat g$,
that $C_2$ is strictly hyperbolic if $\tau_1\ne0$.
(If $\tau_{1}\ne0$, $C_2$ has real eigenvalues
$\lambda_{1}=0,\,\lambda_{2}=-\tau_{1}$, and the left eigenvectors
$\bv_1=(-\tau_1,1)$ and $\bv_2=(-2\tau_1,1)$.
If~$\tau_{1}\equiv0$, i.e., $\calf$ is minimal foliation, the matrix $C_2$ is nilpotent,
hence it is not hyperbolic).
By Corollary~\ref{C-A0}, $\hat g_t= \hat g_0\exp\big(\int_0^t \sigma_2\,dt\big)$,
where $\sigma_2$ exists for $0\le t<T$.

\underline{Remark} that $\tau_1^2-2\tau_2=-(k_1-k_2)^2=const$ along the first family of characteristics $\frac{dx}{dt}=0$, hence the function $k_1-k_2$ does not depend on $t$
(the flow preserves the property ``umbilic free foliation").
Along the second family of characteristics $\frac{dx}{dt}=-\tau_1$ we have
$k_1k_2=(\tau_1^2-\tau_2)/2=const$.
 Notice that (\ref{E-ExtRick}) is reduced to two equations for $k_1$ and $k_2$ with equal RHS's,
\begin{equation*}
  \dt k_i = \dx(k_1 k_2),\qquad i=1,2,
\end{equation*}
hence $\dt(k_2-k_1)=0$, i.e., $k_2-k_1=\psi(x)$ is a known (at $t=0$) function of $x$.
Recall that the system is hyperbolic in the $t$-direction, if $k_1+k_2\ne0$, otherwise, the matrix is nilpotent.
For $k_1$ we have a quasilinear PDE $\,\dt k = \dx(k^2+\psi(x)\,k)$,
whose solution exists for small enough $t$.

\underline{For $n=3$}, the matrix of 3-truncated system (\ref{E-Rictau}) is
$$
 C_3=\left(\begin{matrix}
                     -2\,\tau_{1} & 1        & 0    \\
                     -2\tau_{2} & -\tau_{1}  & \frac43 \\
                    \tau_{1}^3{-}\tau_{3}{-}3\tau_{1}\tau_{2} & \frac32(\tau_{2}{-}\tau_{1}^2) & \tau_{1}\\
           \end{matrix}\right).
$$
Replacing $\tau$-s by $\sigma$-s, see Remark~\ref{R-newton}, we obtain
the characteristic polynomial $P_3=\lambda^{3} + 2\,\sigma_1\lambda^2 + \sigma_1^2\lambda + 4\,\sigma_3$.
Substituting $\lambda=y-\frac23\,{\sigma_1}$ into $P_3$ gives
\begin{equation*}
 P_3={y}^{3}+py+q,\quad
 \mbox{where}\quad p=-(1/3)\,\sigma_{1}^{2}\ \mbox{ and }\ q=4\,\sigma_{3}-({2}/{27})\,\sigma_{1}^{3}.
\end{equation*}
 Depending on the sign of the discriminant $D=(q/2)^2+(p/3)^3$, we have
$$
 D\,\left\{
 \begin{array}{cc}
   > 0, & \mbox{\rm one real and two complex roots}, \\
   < 0, & \mbox{\rm three different real roots},\\
   = 0, & \mbox{\rm one real root with multiplicity three
   in the case } p = q = 0,\\
        & \mbox{\rm or a single and a double real roots when } (\frac13p)^3{=}{-}(\frac12q)^2\ne0.
 \end{array}
 \right.
$$
In our case, $D=\frac4{27}\,\sigma_3(27\,\sigma_3-\sigma_1^3)$.
The condition ``three different real roots"~is
 $D<0 \ \Leftrightarrow \ |\sigma_1|^3>27\,|\sigma_{3}|>0$.\hfill$\square$

\begin{cor}\label{C-tu2}
Let $(M, g_0)$ be a Riemannian manifold,
and $\calf$ a codimen\-sion-1 totally umbilical foliation on $M$ with the normal curvature $\lambda _0$
and a complete unit normal field~$N$. Set at $t=0$

\vskip2mm
 $T=\infty$ if $N(\lambda _0^2)\le 0$ on $M$, and $T=1/[(n-1)\,\sup _{M}N(\lambda _0^2)]$ otherwise.

\vskip2mm\noindent
Then the extrinsic Ricci flow (\ref{E-extRic}) has a unique smooth solution $g_t$ on $M$ for $t\in[0, T)$,
and does not possess one for $t\ge T$.
\end{cor}

\noindent\textbf{Proof}.
The function $\tilde\lambda_t(x)=\lambda(\gamma(x),t)$ along the trajectory $\gamma(x)\ (\gamma(0)=p)$, of $N$ satisfies
(\ref{E-tauAk-L}) with $\psi(\lambda)=4(1-n)\lambda^2$ and
initial value $\tilde\lambda_0(x)=\lambda(\gamma(x),0)$. Then we apply Proposition~\ref{C-tu1}.\hfill$\square$

\vskip2mm
 The Newton transformations of the shape operator can be applied successfully to foliated manifolds, see \cite{rw2}.
By the Cayley-Hamilton theorem, $T_n(b)=0$.
The \textit{extrinsic Newton transformation (ENT)} flow is given~by
\begin{equation}\label{eq1-NT1}
\dt g_t = T_s(b_t).
\end{equation}
In other words, we have $h(b)$ in (\ref{eq1}) with $f_j=(-1)^j\sigma_{s-j}$.
Notice that $\tr T_s(b)=(n-s)\,\sigma_s$ (the trace).
For the ENT flow, (\ref{E-tauAk}) reduces~to
\begin{equation*}
 \dt\tau_i +\frac{i}2\big\{\hskip-1pt N(\sigma_{s})\,\tau_{i-1}
 +\hskip-1pt\sum\limits_{j=1}^{s}\hskip-1pt(-1)^j\big[N(\sigma_{s-j})\,\tau_{i+j-1}
 +\frac{j}{i{+}j{-}1}\,\sigma_{s-j} N(\tau_{i+j-1})\big]\hskip-1pt\big\}=0.
\end{equation*}
From (\ref{eq6}) in what follows with $f_j$ replaced by $(-1)^j\sigma_{s-j}$, we also get
\begin{eqnarray}\label{eq6-NT}
\nonumber
 \dt\sigma_s \eq-\frac{1}{2}\sum\nolimits_{i=0}^{s-1}\sigma_{s-i-1}\big\{ (-1)^{i} N(\sigma_{s})\,\tau_{i}\\
 \plus\sum\nolimits_{j=1}^s (-1)^{i+j}\big[N(\sigma_{s-j})\,\tau_{i+j}+\frac{j}{i+j}\,\sigma_{s-j}\,N(\tau_{i+j})\big]\big\}.
\end{eqnarray}
For $s=1<n$, we have the linear PDE
 $2\,\dt\sigma_1 = -(n-1)\,N(\sigma_1)$
(see also (\ref{E-tauAk}) with $i=1)$ representing the ``unidirectional wave motion"
 $\sigma^t_1(s)=\sigma_1^0(s-t(n-1)/2)$ along any $N$-curve $\gamma(s)$.
For $s=2<n$, (\ref{eq6-NT}) is reduced to the PDE
 $2\,\dt\sigma_2 = [\tau_2-(n-1)\,\tau_1^2]\,N(\tau_1) +\frac n2\,\tau_1\,N(\tau_2)-\frac23\,N(\tau_3)$.

\subsection{EGF with rotational symmetric metrics}

Notice that the EGF preserves rotational symmetric~metrics
 \begin{equation}\label{E-G0ij-6}
 g_t = dx_0^2+\varphi^2_t(x_0)\,ds_{n}^2,\quad
 \mbox{where }\ ds_{n}^2\ \  \mbox{is a metric of curvature 1}.
\end{equation}
The $n$-parallels $\{x_0=c\}$ compose a
(Riemannian) totally umbilical foliation $\calf$ with $N=\partial_0$.
 In this case, $\lambda_t$ can be found from (\ref{E-tauAk-L}), and
 by Proposition~\ref{C-tu1}, the EGF with generating functions $f_j=f_j(\vec\tau)$ can be reduced to
\begin{equation}\label{E-rotsymflow}
 \dt g_t=\psi(\lambda_t)\,\hat g_t,\quad{\rm where}\quad \psi_t(\lambda):=
 \sum\nolimits_{j=0}^{n-1}f_j(n\lambda,n\lambda^2,\ldots,n\lambda^n)\,\lambda^j.
\end{equation}
For simplicity, assume $n=1$, and consider $(M^2,g_t,\calf)$ in biregular foliated coordinates $(x_0,x_1)$. Hence $(g_t)_{00}=1$ and $(g_t)_{11}=\varphi_t^2$.
From (\ref{E-rotsymflow}) we get
\begin{equation}\label{E-G0ij-2mu}
 \dt\varphi_t=(1/2)\,\psi(\lambda_t)\,\varphi_t,
 \ \mbox{ or } \
 \varphi_t=\varphi_0\,e^{(1/2)\int_0^t\psi(\lambda_t)\,dt}.
\end{equation}
By Lemma~\ref{lem-biregG} (for $b_{11}=g_t(A\partial_1,\partial_1)=\lambda_t\,\varphi_t^2$), we also have
 $\lambda_t= -(\varphi_t)_{,0}/\varphi^2_t$.
The Gaussian curvature of $M^2$ is $K_t=-(\varphi_t)_{,00}/\varphi_t$.
 For example, if $\psi(\lambda)=\lambda$, then (\ref{E-tauAk-L}) reduces
to the linear PDE $\dt\lambda+\frac12 N(\lambda)\,{=}\,0$
representing the ``unidirectional wave motion" along any $N$-curve $\gamma(s)$
\begin{equation}\label{eq-unmotion}
 \lambda_t(s)=\lambda_0(s-t/2).
\end{equation}
 In~particular, if $\lambda_0=C\in\RR$, then also
$\lambda_t=C$, and $\varphi_t=\varphi_0\,e^{(t/2)\,\psi(C)}$.

\begin{figure}[ht]
\begin{center}
\includegraphics[scale=.3,angle=0,clip=true,draft=false]{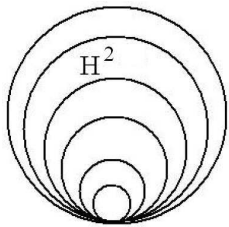}
\caption{\small Foliation of ${\mathbb H}^2$ with constant $\lambda\ne0$.}
\label{F-002}
\end{center}
\end{figure}

Foliations with $\lambda=\mbox{const}\ne0$ exist on a hyperbolic plane. Their leaves are horocycles. On the Poincar\'{e} $2$-disc $B$ the leaves of such foliations are represented by Euclidean circles tangent $\partial B$, Figure~\ref{F-002}.
Orthogonal trajectories to above foliations form foliations by geodesics (i.e., $\lambda=0$).

\begin{ex}\rm
Some of rotational symmetric metrics come from surfaces of revolution in Euclidean 3-space. Evolving them by EGF yields deformations of surfaces of revolution foliated by parallels.
 Revolving a curve $\gamma_t(x)=(X_t(x),\,Y_t(x)$, $Z=0),\ (x\in I)$, around the $x$-axis, we get a surface of revolution $M^2_t\subset\RR^3$,
 $(x,y)\to \big(X_t(x),\ Y_t(x)\cos y,\ Y_t(x)\sin y\big)$.
Assuming that $(X_t'(x))^2+(Y_t'(x))^2=1$ ($x$ is a natural parameter on $\gamma_t$), we obtain the rotational symmetric metric
$g_t = dx^2+Y_t^2(x)\,dy^2$ with $|Y_t'(x)|\le1$.
One may recover $\gamma_t$ from $g_t$ by its curvature $k^2(\gamma_t)=(X_t'')^2+(Y_t'')^2$.

 Revolving the graph of $y=Y(x)$ about the $x$-axis, we get the metric
\begin{equation}\label{eq_g-rev}
 g=(1+{Y'}^{\,2})\,dx^2+Y^2\,dy^2
\end{equation}
on the cone $y^2+z^2=Y^2(x)$.

(i) Revolving a line $\gamma_0(x)=(x\cos\beta,\, x\sin\beta,\, 0)$ about the
$x$-axis, we build the cone
$C_0: y^2+z^2=x^2\tan^2\beta$ (foliated by parallels $\{x=c\}$)
with the first fundamental form
$g_0=dx^2+(x\sin\beta)^2\,dy^2$.
Hence $\varphi_0=x\sin\beta$ and $\lambda_0(x)=-2/x$.
Applying the flow $\dt g_t=\lambda_t\hat g_t$, by (\ref{eq-unmotion}) we obtain
$\lambda_t(x)=-\frac2{x-t/2}$.
The rotational symmetric metric $g_t=dx^2+((x-t/2)\sin\beta)^2\,dy^2$ appears on
the same cone translated across the $x$-axis, $C_t: y^2+z^2=(x-t/2)^2\tan^2\beta$.

(ii) Let us find a curve $Y=Y(X)>0$ such that the metric (\ref{eq_g-rev})
(on the surface of revolution $M^2: [f(Y), Y\cos Z, Y\sin Z]$) has $\lambda_0\,{=}\,\mbox{const}\,{=}\,1$.
Using $\lambda_0=(1/Y(X))\sin\phi$, where $\tan\phi=Y'(X)$,
we get the ODE
 $1=\frac{|Y'(X)|}{Y(X)\sqrt{1+(Y'(X))^2}} \ \Rightarrow \
  \frac{dY}{dX}= \frac{Y}{\sqrt{4+Y^2}}.
 $
\begin{figure}[ht]
\begin{center}
\includegraphics[scale=.4,angle=0,clip=true,draft=false]{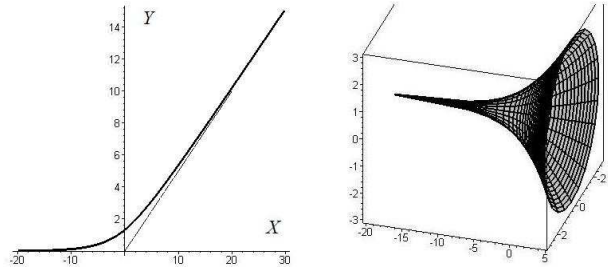}
\caption{\small a) Graph of $X=\log\frac{\sqrt{4+Y^2}-2}{\sqrt{4+Y^2}+2}+\sqrt{4{+}Y^2}$.
\ \ b) Surface of revolution.}
\label{F-001}
\end{center}
\end{figure}
The solution is
 $X=\log\frac{\sqrt{4+Y^2}-2}{\sqrt{4+Y^2}+2}+\sqrt{4+Y^2}+C$, where $C\in\RR$.
The surface $M^2$ looks like a pseudosphere, see Figure~\ref{F-001},
but for $Y\to\infty$ it is asymptotic to the cone $[Y+C,\,Y\cos Z,\,Y\sin Z]$.
The Gaussian curvature of $M^2$ is $K=\frac{-1}{(Y^2+2)^2}<0$, and $\lim\limits_{Y\to\pm\infty}K=0$.
\end{ex}

\subsection{EGF on foliated surfaces}

In this section, $(M^2,g_t)$ is a two-dimensional Riemannian manifold (surface) with a transversally orientable foliation $\calf$ (by curves), $N$ a unit normal to $\calf$, and  $\lambda_t$ the geodesic curvature of the leaves with respect to~$N$.

For $\psi\in C^2(\RR^2)$, the EGF $g_t$ of type (b) on $(M^2,\calf)$ is a solution to the PDE $\dt g_t=\psi(\lambda_t,t)\,\hat g_t$,
where $\lambda_t$ obeys the PDE $\dt\lambda_t+\frac12\psi'_\lambda(\lambda_t,t)N(\lambda_t)=0$ (see (\ref{E-tauAk}) for $i=0$).
By Corollary~\ref{C-A0}, if $\psi'(\lambda_0)\ne0$, see the condition $(H_1)$ with $n=1$, then there is a unique local smooth solution $\lambda_t$ for $0\le t<T$ with initial value $\lambda_0$ determined by $g_0$,
and $\hat g_t= \hat g_0\exp(\int_0^t \psi'_\lambda(\lambda_t,t)\,dt)$ holds.
\newline
For compact $M^2$, we have $\dt(d\vol_t)= \frac12\tr h(b_t) d\vol_t$, see \cite{to}.
Hence, the volume $\vol_t:=\int_M d\vol_t$ of $g_t$ satisfies to
\begin{equation}\label{E-intpsi}
 \dt\vol_t=\frac12\int_M \psi(\lambda_t,t)\,d\vol_t.
\end{equation}

In order to estimate the time interval due to Proposition~\ref{C-tu1}, suppose that the function $\psi\in C^2(\RR)$ does not depend on $t$. Define $T=\infty$ if $N(\psi(\lambda_0))\ge0$ on $M$, and $T=-2/\inf_{M}N(\psi(\lambda_0))$ otherwise.
 Then the EGF $\dt g_t=\psi(\lambda_t)\,\hat g_t$ has a unique smooth solution $g_t$ on $M^2$ for $t\in[0, T)$, and does not possess one for $t\ge T$.
 If, in addition to above $(H_1)$, $\psi'(\lambda_0)N(\lambda_0)\ge0$ holds, then the solution exists for all $t\ge0$ ($T=\infty$).

\begin{prop}\label{L-nablaNNt}
The Gaussian curvature of the EGF of type (b) on $(M^2,\calf)$ is given by the formula
\begin{equation}
\label{E-gaussM2}
 K_t = \Div(e^{-\int_0^t\psi(\lambda_t,t)\,dt} \nabla^0_NN) +N(\lambda_t) -\lambda_t^2.
\end{equation}
\end{prop}

\noindent\textbf{Proof}.
Define a self-adjoint $(1,1)$-tensor $h(A):=\sum\nolimits_{j=0}^{n-1} f_j\,A^ j$ dual to $h$ of (\ref{eq1}).
We use (\ref{eq2}) to compute for any $X\in T\calf$
\begin{eqnarray*}
 g(\dt(\nabla^t_NN), X) \eq \frac12\big[2(\nabla^t_N h_t)(X,N)-(\nabla^t_X h_t)(N,N)\big]\\
 =-h(\nabla^t_NN, X)\eq -\sum\nolimits_{j=0}^{n-1} f_j\,g(A^ j(\nabla^t_NN), X)
 =-h(A)(\nabla^t_NN).
\end{eqnarray*}
Hence we have the general relation
  $\dt(\nabla^t_NN) = -h(A)(\nabla^t_NN)$,
 which for $n=1$ looks~as
  $\dt(\nabla^t_NN)\,{=}-\psi(\lambda_t,t)\nabla^t_NN$.
Integrating above~yields
 \begin{equation}
\label{E-nablaNNt}
  \nabla^t_NN = e^{-\int_0^t\psi(\lambda_t,t)\,dt} \nabla^0_NN.
 \end{equation}
The formula (\ref{E-gaussM2}) for $K=\Ric(N,N)$ is a consequence of (\ref{E-nablaNNt}) and the formula for ${\rm div}(\nabla_N N)+{\rm div}(H)$ (of a codimension-one foliation), see \cite{wa1},
\begin{equation*}
 {\rm div}(\nabla_N N) =\Ric(N,N) + \tau_2 -N(\tau_1).\qquad\square
\end{equation*}

\hskip-3mm
(a) Let $\psi=1$. The solution to $\dt g_t=\hat g_t$ is $\hat g_t=e^t\hat g_0\ (t\ge0)$. From $\dt\lambda=0$ we get $\lambda_t=\lambda_0$.
 By (\ref{E-gaussM2}), the Gaussian curvature is
 $ K = e^{-t}\Div(\nabla^0_NN)+N(\lambda_0) -\lambda_0^2$.
 There is limit $K_\infty=N(\lambda_0)-\lambda_0^2$.

\vskip.5mm\hskip-3mm
(b) Let $\psi=\lambda$, i.e., $\dt g=\hat b_1$. Then
$\lambda_t(s)=\lambda_0(s+\frac t2)$ along any $N$-curve $\gamma(s)$
in the $(t,s)$-plane.
From 
$\dt g_t=\lambda_0(s+\frac t2)\hat g_t$ we get
$\hat g_t(s)\,{=}\,\hat g_0(s)e^{\,\int_0^t \lambda_0(s+\xi/2)\,d\xi}$ $(t\in\RR)$.
For compact $M^2$, by (\ref{E-intpsi}) we have $\vol_t=\mbox{const}$.

\vskip.5mm\hskip-3mm
(c) Let $\psi=\lambda^2$. Then $\dt\lambda+\lambda\,N(\lambda)=0$ (the inviscid Burgers equation). If $N(\lambda_0)\ge0$ (for $t=0$), then the solution $\lambda_t$ exists for all $t\ge0$.

\begin{ex}\label{ex-Reeb1}\rm
Let a function $f\in C^2(-1,1)$ has vertical asymptotes $x=\pm1$.
Consider the foliation $\calf$ in the closed strip, whose leaves
are $L_{\pm}=\{x=\pm1\}$ and
 $L_s(x)= \{(x,\ f(x)+s),\ |x|<1\}$, where $s\in\RR$.
The normal $N$ at the origin is directed along $y$-axis.
The tangent and normal to $\calf$ unit vector fields (on the whole strip) are
$X=[\cos\alpha(x),\,\sin\alpha(x)],\ N=[-\sin\alpha(x),\,\cos\alpha(x)],$
where $\alpha(x)$ is the angle between the leaves $L_s$ and the $x$-axis
at the intersection points.
Indeed, $f$ and $\alpha$ are related~by
\begin{equation*}
 f'(x)=\tan\alpha(x)\quad\mbox{and}\quad
 \cos\alpha=[1+(f')^{2}]^{-1/2},\ \ \sin\alpha=f'\,[1+(f')^{2}]^{-1/2}.
\end{equation*}
The curvature of $L_s$ is
  $\lambda_0(x)={f''(x)}[1+(f'(x))^{2}]^{-3/2}=\alpha'(x)\cdot|\cos\alpha(x)|$,  where $|x|<1$.
The $N$-curves through the critical points of $f$ are vertical,
and divide $\Pi$ into sub-strips.
Typical foliations in the strip $|x|<1$ with one vertical trajectory $x=0$ are the following two:

(i) $f$ has exactly one strong minimum at $x=0$.

(ii) $f$ is monotone increasing with one critical point $x=0$;

\noindent
Taking $f=\frac1{10}[e^{x^2/(1-x^2)}-1]$ or $\alpha(x)=\frac\pi2 x$
for (i), we get the Reeb foliation.

\noindent
For (ii) one may take $f=\tan(\frac\pi2 x)$, or $\alpha(x)=\frac\pi2 x^2$.

\vskip1.5mm
Let $\psi=\psi(\lambda_t)$, where $\lambda_t(x)$ is known for a positive time interval $[0,\eps)$, see Proposition~\ref{C-tu1}.
 We use $X\in T\calf$ and normal $N$ to represent the standard frame
  $e_1=\cos\alpha(x)X-\sin\alpha(x)N$, $e_2=\sin\alpha(x)X+\cos\alpha(x)N$
in the $(x,y)$-plane.
 By $\hat g_t=\hat g_0\,e^{\int_0^t\psi(\lambda_t(x))\,dt}$, we get
 $g_t(X,X)=e^{\int_0^t\psi(\lambda_t(x))\,dt},\
 g_t(X, N)=0$, and
 $g_t(N, N)=1$.
The $g_t$-scalar products of the frame $\{e_1, e_2\}$ are
\begin{eqnarray*}
 E_t\eq g_t(e_1,e_1)=\sin^2\alpha+\cos^2\alpha\,e^{\int_0^t\psi(\lambda_t(x))\,dt},\\
 F_t\eq g_t(e_1,e_2)=\sin\alpha\cos\alpha[e^{\int_0^t\psi(\lambda_t(x))\,dt}-1],\\
 G_t\eq g_t(e_2,e_2)=\cos^2\alpha+\sin^2\alpha\,e^{\int_0^t\psi(\lambda_t(x))\,dt}.
\end{eqnarray*}
The Gaussian curvature $K_t$ of the metric $g_t=E\,dx^2+2F\,dx\,dy+G\,dy^2$ is 
\begin{equation*}
\begin{array}{ccc}
 K_t=-\frac{1}{2\sqrt{E G-F^2}}\,
 \Big[\dx\big(\frac{\dx G-\dy F}{\sqrt{E G-F^2}}\big)
 +\dy\big(\frac{\dy E-\dx F}{\sqrt{E G-F^2}}\big)
 -\frac{1}{4(E G-F^2)^2}
 \left|\begin{smallmatrix}
     E & \dx E & \dy E \\
     F & \dx F & \dy F \\
     G & \dx G & \dy G \\
   \end{smallmatrix}\right|\Big]
\end{array}
\end{equation*}
which in our case, when the coefficients of $g_t$ do not depend on $y$-coordinate, reads as
 $K_t=\frac{-1}{2\sqrt{E_t G_t-F_t^2}}\,\dx\big(\frac{\dx G_t}{\sqrt{E_t G_t-F_t^2}}\big)$.
We have $E_t G_t-F_t^2=e^{\int_0^t\psi(\lambda_t(x))\,dt}$.

The $N$-curves satisfy ODEs
$\,dx/dt=-\sin\alpha(x),\ dy/dt=\cos\alpha(x)$.
From the first of above ODEs for $N$-curves, we deduce
the implicit formula $t=-\hskip-1pt\int_x^{\phi_t(x)}\hskip-1pt\frac{dx}{\sin\alpha(x)}$
for local diffeomorphisms $\phi_t(x)\ (|x|<1,\,t\ge0)$.

\underline{Suppose that $\psi(\lambda)=\lambda$}.
Since $\lambda_t(s)=\lambda_0(s+t/2)$ is a simple wave along $N$-curves, we have
$\lambda_t(x)=\lambda_0(\phi_{t/2}(x))$.
For example, $\lambda_t(0)=\lambda_0(0)$ for all $t\ge0$.
Substitution $E_t,F_t$ and $G_t$ into above formula for $K_t$ yields
\begin{equation}\label{E-Kiii}
 \hskip-7mm
 \begin{array}{ccc}
 && K_t=\frac18(\cos(2\alpha)-1)
 \big[2\int_{0}^{t}\!\frac{\partial^{2}}{\partial{x}^{2}}\lambda_t\,dt
 +\big(\int_{0}^{t}\!{\frac{\partial}{\partial x}}\lambda_t\,dt\big)^{2}\big]
 -\big[\cos(2\alpha)(\alpha')^{2}\\
 &&+\frac12\,\sin(2\alpha)\,\alpha''\big]
 \big[1-{\rm e}^{-\int_{0}^{t}\!\lambda_t\,dt}\big]
 -\frac14\,\sin(2\alpha)\,\alpha'
 \int_{0}^{t}\!{\frac{\partial}{\partial x}}\lambda_t\,dt
 \big[3+{\rm e}^{-\int_{0}^{t}\!\lambda_t\,dt}\big].
\end{array}
\end{equation}
Since $\alpha(0)=0$ and $\lambda_t(0)=\lambda_0(0)$, one has
$K_t(0)=-(\alpha'(0))^{2}\big[1-{\rm e}^{-t\,\lambda_0(0)}\big]$.
\newline
\underline{For (i)}, we get $\alpha'(0)>0$ and $\lambda_0(0)>0$, so $K_t(0)<0$ for $t>0$. Since $\lim\limits_{t\to\infty}\phi_t(x)=0$ for $|x|\le 1$, there also exists $\lim\limits_{t\to\infty}\lambda_t(x)
 =\lambda_0(0)>0$. Hence for any $x\in[-1,1]$ there is $t_x>0$ such that
$K_t(x)<0$ for $t>t_x$.

\underline{For (ii)} with $\alpha(x)=(\pi/2)x^2$, we have $\alpha'(0)=0$ and $\alpha''(0)=\pi$.
Moreover, $\lambda_0(0)=0$ and $\lambda'_0(0)=\pi\ne0$.
Since $\lim\limits_{t\to\infty}\phi_t(x)=0$ for all $0<x\le 1$, there also exists $\lim\limits_{t\to\infty}\lambda_t(x)=\lim\limits_{t\to\infty}\lambda_0(\phi_{t/2}(x))
=\lambda_0(0)=0$.
By (\ref{E-Kiii}) we have $K_t(0)\equiv0$ and the series expansion
  $K_t(x)=-\frac52\,\pi^{2}\,\int_0^t\frac{\partial}{\partial x}\lambda_t\,dt\,x^3
  +O_t\big({x}^{4}\big)$.
We conclude that there exists $t_0>0$ such that $K_t(x)$ for $t > t_0$ changes its sign when we cross the line $x=0$.
\end{ex}

\begin{ex}\rm
 Consider a foliation $\calf$ by circles $L_\rho=\{\rho=c\}$ in the ring
 $\Omega=\{c_1\le\rho\le c_2\}$ for some $c_2>c_1>0$
 with polar coordinates $(\rho,\theta)$.
 Then $X=\partial_\theta$ and $N=\partial_\rho$ are tangent and normal vector fields to the foliation. The metric is $ds^2=d\rho^2+G_t(\rho)\,d\theta^2$.
 Notice that $\lambda_0(\rho)=1/\rho$.
 Since $\dt G_t= \psi(\lambda_t(\rho),t)\,G_t$, we have
 $G_t=\rho^2 \exp(\int_0^t \psi(\lambda_t(\rho),t)\,dt)$.
 From the formula for Gaussian curvature of Example~\ref{ex-Reeb1} we have
 $K_t=-\frac{1}{2\sqrt{G_t}}\,\partial_\rho\big(\frac{\partial_\rho G_t}{\sqrt{G_t}}\big)$.
 Let \underline{$\psi=\lambda$}. Then $\lambda_t(s)= \lambda_0(s+t/2)$ on the $N$-curves.
 For the foliation by circles $\rho=c$ we have $\lambda_t(\rho)= (\rho+t/2)^{-1}$.
 The Gaussian curvature $K_t\equiv0$.
\end{ex}

\end{document}